\documentclass[11pt,letterpaper]{amsart}
\usepackage{latexsym,amsfonts,amssymb,amsmath,amsthm,amscd}
\newtheorem{theorem}{Theorem}[section]
\newtheorem{lemma}[theorem]{Lemma}

\newtheorem{corollary}[theorem]{Corollary}

\theoremstyle{definition}
\newtheorem{example}[theorem]{Example}

\newtheorem{definition}[theorem]{Definition}

\newtheorem{remark}[theorem]{Remark}
\usepackage[utf8]{inputenc}
\usepackage[T1]{fontenc}
\begin{document}

\begin{center}
{\bf{Macdonald-Hurwitz Number}}\\
Quan Zheng \footnote{E-mail: quanzheng2100@163.com}\\
Mathematics College, Sichuan University, 610064, \\
 Chengdu, Sichuan, PRC\\
 November 11, 2022
\end{center}

\begin{center}
\vskip 0.4in
{\bf Abstract}\\[\baselineskip]

\parbox{14cm}{\quad
 Inspired by J. Novak's works on the asymptotic behavior of the BGW and the HCIZ matrix integrals  \cite{[N0]} and by the algebraic and geometric properties of the Hurwitz numbers \cite{[IP]}, \cite{[LZZ]}, \cite{[LR]}, \cite{[OP]}, \cite{[Z1]}, and by the symplectic surgery theory of the relative GW-invariants \cite{[IP]}, \cite{[LR]}, using the elements of the transform matrix from the integral Macdonald function with two parameters to the homogeneous symmetric power sum functions  \cite{[M]}, we have constructed the Macdonald-Hurwitz numbers.  As an application, we have constructed a series of new genus-expanded cut-and-join differential operators, which can be thought of as the generalization of the Laplace-Beltrami operators and have the genus-expanded integral Macdonald functions as their common eigenfunctions. We have also obtained some generating wave functions of the same degree, which are generated by the Macdonald-Hurwitz numbers and can be expressed in terms of the new cut-and-join differential operators and the initial values. Another application is that we have constructed a new commutative associative algebra  $(C(\mathbb{F}[S_{d}]),\circ_{q,t})$ (referring to the last section (6)).
  By taking the limit along a special path $\eta(A|B)$ (referring to the formulas (\ref{p1}), (\ref{p2})), we specialize $(C(\mathbb{F}[S_{d}]),\circ_{q,t})$ to be a   commutative associative algebra  $(C(\hat{\mathbb{F}}[S_{d}]),\circ_{A|B})$, which will be proven to be isomorphic to
  the middle-dimensional $\mathbb{\mathbb{C}^*}$-equivalent cohomological rings via the Jack functions over the Hilbert scheme points of $\mathbb{C}^2$ constructed by W. Li, Z. Qin, and W. Wang in \cite{[LQW2]}.}

\end{center}
\indexspace
\renewcommand\contentsname{\Large Contents}
\tableofcontents

\section {Introduction}\label{mi}
Our primitive motivation for the paper is to understand and generalize J. Novak's works about the asymptotic behavior of the BGW and the HCIZ matrix integrals \cite{[N0]}. To do this, our first step is to set up a link between the integral Macdonald functions with the theory of the Riemann surfaces by applying the Macdonald-Hurwitz numbers defined in the paper, referring to the definition (\ref{definition}). In our views, the asymptotic expansions of the BGW and the  HCIZ matrix integrals are actually related to the theory of the Riemann surface, because these integrals can be expressed by the Hurwitz number \cite{[H]} or the Schur functions, both of which have the orthogonality \cite{[M]}, \cite{[Z1]}. More generally, if an integral can be expressed by the Jack functions with a Jack parameter $\alpha $ (including the zonal functions $\alpha=2$, the Schur functions $\alpha=1,$ and the symplectic zonal functions $\alpha=\frac{1}{2}$),  or the (integral) Macdonald functions with two parameters  $q, t,$ \cite{[M]}, then the asymptotic expansions should be actually related to the theory of the Riemann surfaces due to the Macdonald-Hurwitz numbers,  or to the orthogonality of those functions,  seeing  \cite{[M]} or the orthogonal/cutting lemma (\ref{orth}) in the paper. We plan to further investigate this in our next step.

If we say that the differential operators, eigenfunctions, and eigenvalues are the three musketeers playing important roles in mathematics and physics, then the classical cut-and-join operators \cite{[GJ1]}, \cite{[GJ2]}, \cite{[GJ3]}  and the Schur functions \cite{[M]} are the pair of rangers who are playing a stage drama of love and hate. They are flashing from time to time in the classical \cite{[H]} or all kinds of refined (such as \cite{[CD]}, \cite{[DKPS]}) Hurwitz Enumeration Problem and the integrable systems, the Hodge integral,  the harmonic analysis in the Lie group,  etc  \cite{[GJ1]}, \cite{[GJ2]}, \cite{[GJ3]}, \cite{[LZZ]}, \cite{[OP]}. Meanwhile, the Hurwitz theory serves as a great matchmaker to connect the theory of the Riemann surface with the theory of the BGW and the HCIZ integrals \cite{[HC]}, \cite{[IZ]}. The people pursue the topology/genus expansion property of a differential operator or a function, which is demonstrated clearly and vividly in the study of the limit of the BGW integrals and the HCIZ integrals \cite{[N0]}, \cite{[N]}. We will continue to adhere to this principle in this paper.

 Inspired by the algebraic expressions of the Hurwitz numbers, which are mainly expressed by the elements of the transform matrix from the Schur functions $s_\lambda$ to the homogeneous symmetric power sum functions $p_\lambda$, and by the symplectic surgery theory of the relative GW-invariants developed by A.M. Li and Y.B. Ruan \cite{[LR]}, E. Ionel and T. Parker  \cite{[IP]}, using the elements of transform matrix from the integral Macdonald functions $J_\lambda(x; q,t)$ with two parameters $(q,t)$ \cite{[M]} to the homogeneous symmetric power sum functions $p_\lambda$, we have constructed the Macdonald-Hurwitz numbers. As an application, we obtain the new genus-expanded cut-and-join operators (\ref{cut}) similar to \cite{[AMMN]}, \cite{[MMN]}, \cite{[Z1]},\cite{[Z2]}, which have the genus-expanded integral Macdonald functions (\ref{11}) as the common eigenfunctions.  The genus-expanded cut-and-join operators can be thought of as the generalization of the Laplace-Beltrami operators. \\

As a special case, let us see what has happened to the zonal polynomials \cite{[J64]}, \cite{[J]}, \cite{[S]}.
From the view of the geometry, it is well known that the zonal polynomials are   eigenfunctions  of the  Laplace-Beltrami operator   \cite{[Hel]}, \cite{[J]}
\begin{equation}
\Delta=\frac{1}{det(g_{ij})}\sum_{i,j=1}^N \frac{\partial}{\partial x_i} \sqrt{det(g_{ij})}g^{ij} \frac{\partial}{\partial x_j},
\end{equation}
where $x_1,x_2,\cdots,x_N$  are coordinates of a point in a space $M$ with Riemannian metric differential form
\begin{eqnarray}
&& (ds)^2=\sum_{i,j=1}^{N}g_{ij}dx_i dx_j\\
&& (g^{ij})=(g_{ij})^{-1}
\end{eqnarray}
Especially, let $GL(N,{\mathbb{R}}), O(N,{\mathbb{R}})$ be the real general linear group and the real orthogonal linear group, then we consider the homogenous space $M \cong GL(N,{\mathbb{R}})/O(N,{\mathbb{R}})$, which is isomorphic to the space of $N \times N$ real positive definite symmetric matrices $X$ (by sending a nonsingular matrix  $V$  to $VV^{\prime}$)  with congruence invariant Riemannian metric  \cite{[J]}, \cite{[Hans Maass]}
\begin{eqnarray}
(ds)^2=tr(X^{-1}dX X^{-1}dX).
\end{eqnarray}
Then the Laplace-Beltrami operator concerned with   the latent root $x_1,\cdots,x_N$ of $X$ is the following  \cite{[J]}
\begin{equation}\label{Delta}
\hat{D}:=\sum_{i=1}^N x_i^2\frac{\partial^2}{\partial x_i^2} +
\sum_{i=1}^{N}\sum_{j=1,j\neq i}^{N}x_i^2(x_i-x_j)^{-1}\frac{\partial}{\partial x_j}-\sum_{i=1}^{N}\frac{1}{2}(N-3)x_i\frac{\partial}{\partial x_i}
\end{equation}
where the last term will not change a homogeneous polynomial as eigenfunction but the eigenvalue, hence we leave this term out as \cite{[C]}, \cite{[CD]}, \cite{[D]}
\begin{equation}\label{laplace2}
D:=\sum_{i}^N x_i^2\frac{\partial^2}{\partial x_i^2} +
\sum_{i=1}^{N}\sum_{j=1,j\neq i}^{N}x_i^2(x_i-x_j)^{-1}\frac{\partial}{\partial x_j}
\end{equation}

In view of physics, the zonal polynomials are the wave functions and the eigenvalues are the energy eigenvalues \cite{[C]}.  To deal with the $Schr\ddot{o}dinger$ equations for the one-dimensional N-body problem, the physicist  F. Calogero \cite{[C]} proposed a similar Laplace-Beltrami operator as the formula (\ref{Delta}) in 1969.

Let us introduce the time-variables or Miwa variables $p_1,p_2,\cdots, $ (the power sum of the laten roots $x_1,\cdots,x_N$ of $X$), i.e., for any positive integer $i\geq 1$,
\begin{equation}
p_{i} = \sum_{j=1}^{N}x_j^{i},
\end{equation}
then the Laplace-Beltrami operator (\ref{laplace2}) acting on  the functions in the time-variables such  as $f(p_1,p_2,\cdots)$ is equivalent to the classical  cut-and-join operator with partition type $(1^{d-2}2)$ \cite{[GJ1]},\cite{[LZZ]},\cite{[MMN]} as
\begin{eqnarray}\label{opp}
D=\sum_{k,l\geq 1}\left[klp_{k+l}\frac{\partial^2}{\partial p_k \partial p_l} +(k+l)p_kp_l\frac{\partial}{\partial p_{k+l}}\right]
\end{eqnarray}
up a multiple $\frac{1}{2}$, which will again not change the eigenfunctions but eigenvalues.

Meanwhile, the zonal polynomial is closely concerned with the irreducible character of the general real linear group $GL(N, \mathbb{R})$,  \cite{[J60]}, \cite{[M]}.

In this paper, we have extended the above theory to the genus-expanded cut-and-join operators and the genus-expanded integral Macdonald functions by defining the Macdonald-Hurwitz numbers, such that the latter is the common eigenfunctions of the former.
Let $d\geq 0$ be non-negative integer, and $\lambda, \Delta$  be the partition of $d,$ which is denoted by $\lambda, \Delta \vdash d.$
Let $J_\lambda(x;q,t)$ be the integral Macdonald functions with two parameters $q,t,$ \cite{[M]}, which can be expressed by  the power sum functions $p_\Delta$,
say, \cite{[M]}
\begin{equation}
J_{\lambda}(x; q, t)=\sum_{\Delta \vdash d}a_{\lambda}(\Delta)(q,t)p_{\Delta},
\end{equation}
where
\begin{equation}
p_\Delta=\prod_{i}p_{\delta_i}
\end{equation}
if $\Delta=(\delta_1,\delta_2, \cdots),$ and $a_{\lambda}(\Delta)(q,t) \in {\mathbb{ F}}:=\mathbb{Q}(q,t),$ the field of the rational function with the rational coefficients in $q, t,$ especially, we denote $dim\lambda:=a_\lambda(1^{|\lambda|})(q,t).$

Let $\Sigma^h$ be a closed connected Riemann surface of genus $h\geq 0,$ and $\Sigma^g$
 a closed connected Riemann surface of genus $g \geq 0$.
 For a general given points set $\{q_1,\cdots,q_k\}\in \Sigma^h$, which is called the set of
 branch points, we call a
holomorphic map
$f:\Sigma^g\to \Sigma^h$
a ramified
covering of $\Sigma^h$ of degree $d\ge 0$ by $\Sigma^g$ with a ramification type $(\Delta_{1},\cdots,\Delta_{k})$, if
the preimages of $f^{-1}(q_i)=\{p_i^{1},\cdots,p_i^{l_i}\}$  with orders $\Delta_i = (\delta^i_1,\cdots,\delta^i_{l_i})\vdash d$ for $i=1,\cdots,k,$  respectively, where $l_i=l(\Delta_i)$ is the length of the partition $\Delta_i.$
{\definition  (Definition \ref{definition})
We  define the connected genus $g \geq 0$  Macdonald-Hurwitz number with degree $d > 0$ over the  closed  connected  Riemann surface  $\Sigma^h$  with the ramification type $(\Delta_{1},\cdots,\Delta_{k})$ as
\begin{eqnarray}
 &&MH_{h}^{g,d}(\Delta_{1},\cdots,\Delta_{k})=\\
&&\sum_{\lambda \vdash d}
(\frac{dim^2\lambda}{j_{\lambda}(q,t)})^{1-h}\frac{a_{\lambda}(\Delta_1)(q,t)}{dim \lambda}\cdots \frac{a_{\lambda}(\Delta_k)(q,t)}{dim \lambda},\nonumber
\end{eqnarray}
under the condition
\begin{equation}
(2-2h)d-(2-2g)=\sum_{i=1}^{k}(d-l(\Delta_i)),
\end{equation}
otherwise, we define
\begin{eqnarray}
MH_{h}^{g,d}(\Delta_{1},\cdots,\Delta_{k})=0,
\end{eqnarray}
and if we take $d=0,$   it is convenient to  define
\begin{eqnarray}
MH_{h}^{g,0}(\Delta_{1},\cdots,\Delta_{k})=0.
\end{eqnarray}
where  $j_\lambda(q,t)$ is defined by the formula (\ref{j}).
}

Taking the rule of the exponential-logarithmic correspondence, we can define the disconnected Macdonald-Hurwitz numbers, which are obtained by taking the product of all connected the Macdonald-Hurwitz numbers module their automorphisms, and then summing up all possibilities, as described in section (3).

We can  iteratively  calculate the Macdonald-Hurwitz numbers  by decreasing the genus of the target surface and the numbers of the branch points,
{\theorem\  (Theorem\ref{decrease})
\begin{itemize}
\item[(1)] Decrease the genus of the target surface directly, i.e., if $h\geq 1,$ we have
\begin{eqnarray}{}
MH_{h}^{g,d}(\Delta_{1},\cdots,\Delta_{k})=\sum_{\Delta \vdash d}MH_{h-1}^{g,d}(\Delta_{1},\cdots,\Delta_{k},\Delta,\Delta)z_{\Delta}(q,t),
\end{eqnarray}
where $z_{\Delta}(q,t)$ is defined by the formula (\ref{43});
\item[(2)] (Associativity) Decrease the genus of the target surface and the numbers of the branched points by cutting it into two parts, i.e, for $h_1+h_2=h, h_1\geq 0, h_2 \geq 0,$ and  any positive integer $k> l\ge 1$, we have
\begin{eqnarray}
&& MH_{h}^{g,d}(\Delta_{1},\cdots,\Delta_{k})\\
&=&
\sum_{\Delta}MH_{h_1}^{g_1,d}(\Delta_{1},\cdots,\Delta_{l},\Delta)z_{\Delta}(q,t) \nonumber\\
 &&\times MH_{h_2}^{g_2,d}(\Delta,\Delta_{l+1},\cdots,\Delta_{k}),\nonumber
\end{eqnarray}
\end{itemize}
where the genus $g,g_1,g_2$ are determined by the formula $(\ref{hurwitz})$.}

{\remark{
We have noticed the following interesting facts
(see Corollary (\ref{example0}) and Remark (\ref{oo}) and the formula (\ref{43})):
\begin{equation}\label{l1}
MH_{0}^{0,1}((1),(1))=\frac{1}{z_{(1)}(q,t)}=\frac{1-t}{1-q},
\end{equation}
\begin{equation}\label{l2}
MH_{0}^{0,1}((1))=\frac{1-t}{1-q},
\end{equation}
thus we can determine $\frac{1-t}{1-q}$-curves passing through two points or one point.
 }}

Let $\hbar$ be indeterminate, we have

{\definition (Definition \ref{opd}) {Following the idea of \cite{[Z1]}, for any  partition $\Delta \vdash d$,  we define   a genus-expanded cut-and-join operator $
D(\Delta,\hbar)$ (referring to the formula (\ref{cut}))
\begin{equation}
D(\Delta,\hbar)=\sum_{\Gamma,\Gamma^{\prime}\vdash d}
z_{\Gamma^\prime}(q,t)\hbar^{d+l(\Gamma^{\prime})-l(\Delta)
-l(\Gamma)}MH_0^{g,d}(\Gamma^{\prime},\Delta,\Gamma)
p_{\Gamma}\frac{\partial}{\partial p_{\Gamma^{\prime}}},
\end{equation}
where  genus $g$ satisfies:
\begin{equation}
2d-(2-2g)=(d-l(\Gamma^{\prime}))+(d-l(\Delta))+(d-l(\Gamma)),
\end{equation}}}
and $\frac{\partial}{\partial p_{\Gamma^{\prime}}}$ is defined by the formula (\ref{partial}).

Let $J_\lambda(x;q,t,\hbar)$ be the genus-expanded integral Macdonald functions (referring to the formula (\ref{11})), we have

{\theorem (Theorem \ref{eigen}) For any partition $\Delta, \lambda \vdash d,$ the following equality holds:
 \begin{equation}
 D(\Delta,\hbar)J_\lambda(x;q,t,\hbar) = \hbar^{d-l(\Delta)}a_\lambda(\Delta)(q, t) J_\lambda(x; q, t, \hbar),
\end{equation}
i.e., the genus-expanded integral Macdonald functions  $J_\lambda(x; q,t,\hbar)$
 are  the common  eigenfunctions of  the genus-expanded cut-and-join operators $D(\Delta,\hbar)$
and $\hbar^{d-l(\Delta)}a_\lambda(\Delta)(q, t)$  are the corresponding eigenvalues.}

For a given $d$, the genus-expanded cut-and-join operators, which act on  the functions in the time-variables $p = (p_1, p_2,\cdots )$, are closed under the	multiplication, i.e.,

{\theorem (Theorem \ref{operator})For any partitions $\Delta_1,\Delta_2 \vdash d,$ we have
\begin{equation}
 {D(\Delta_1,\hbar) D(\Delta_2,\hbar) = \sum_{\Delta_3 \vdash d}\hbar^{d-l(\Delta_1)-l(\Delta_2)+l(\Delta_3)}
C_{\Delta_1\Delta_2}^{\Delta_3} D(\Delta_3,\hbar) },
\end{equation}
where
\begin{equation}
C_{\Delta_1\Delta_2}^{\Delta_3}=z_{\Delta_3}(q,t)MH_{d}^{g}(\Delta_{1},\Delta_{2},\Delta_{3})
\end{equation}
is a new structure constant
of the center subalgebra $C(\mathbb{F}[S_{d}])$  of the group algebra $\mathbb{F}[S_{d}]$ defined by the formula (\ref{constant}).}\\

\remark{The terminology "new" in the paper is relative to the classical Hurwitz theory \cite {[MMN]}, \cite {[OP]}, \cite{[Z1]}, etc.}

Let $\Phi_h\{\hbar|(u_1,\Delta_1),\cdots,(u_n,\Delta_n)|p^{(1)},\cdots,p^{(k)},p\}$ be the generating wave function (referring to the formula (\ref{99})), then for any $i\in \{1,\cdots,n\}, $ we have

{\theorem (Theorem \ref{de})
\begin{eqnarray}
&&\nonumber \frac{\partial\Phi_h\{\hbar|(u_1,\Delta_1),\cdots,(u_n,\Delta_n)|p^{(1)},\cdots,p^{(k)},p\}}{\partial u_i} \nonumber\\
& =&D(\Delta_i,\hbar)\Phi_h\{\hbar|(u_1,\Delta_1),\cdots,(u_n,\Delta_n)|p^{(1)},\cdots,p^{(k)},p\}.
\end{eqnarray}}

Furthermore, let $S_d$ be the symmetric group of the degree $d,$  then we have applied the Macdonald-Hurwitz numbers to determine  a $\circ_{q,t}$-product on the central subalgebra  $C(\mathbb{F}[S_d])$ of the  group algebra $\mathbb{F}[S_d]$, which can be thought of as a $(q,t)$-deformation product of $C(\mathbb{C}[S_d])$ if we extend the coefficient field $\mathbb{C}$ to the rational functions fields  $\mathbb{F}:=\mathbb{Q}(q,t)$, as described in the last section (6),i.e. we have

{\theorem (Theorem \ref{ca}){ For any partition $\Delta \vdash d, $ let $C_{\Delta} \in C(\mathbb{C}[S_{d}])$ be the classical central element of $C(\mathbb{C}[S_{d}])$, then  $\{C_{\Delta}| \Delta \in{\mathbb{P}}_d\}$ also forms a ${\mathbb{F}}$-basis of $C(\mathbb{F}[S_{d}])$.  By linearly extending,  $C(\mathbb{F}[S_{d}])$ is a new commutative associative algebra with the unit $C_{(1^d)}$ under the following product,
\begin{equation}
C_{\Delta_1}\circ_{q,t}C_{\Delta_2}= C_{\Delta_1\Delta_2}^{\Delta_3}C_{\Delta_3},
\end{equation}}
where
\begin{equation}
C_{\Delta_1\Delta_2}^{\Delta_3}=z_{\Delta_3}(q,t)MH_{d}^{g}(\Delta_{1},\Delta_{2},\Delta_{3})
\end{equation}
is the new structure constant.}

{\remark{ In a series of the papers  \cite{[BG]}, \cite{[BP]}, \cite{[OP2]}, \cite{[OP3]}, \cite{[PT]}, the authors of those papers have constructed the local Gromov-Witten invariants of the curves,  Donaldson-Thomas invariants, the quantum cohomology of the rank 2 bundles over the Riemann surfaces or the Hilbert scheme points of $\mathbb{C}^2$,  the orbifold quantum cohomology of the symmetric product of $\mathbb{C}^2$,  the equivariant version of those invariants, which are proven to be equivalent and closely related with the Jack functions and the  Macdonald functions, it is an interesting problem to research the relationships between those invariants with our Macdonald-Hurwitz numbers, and we hope to come back to pursue the problem in the future. The author would like to thank an anonymous expert for his communications about the above information.  However,   by taking the limit along a special path $\eta(A|B)$ (the formula (\ref{p1}), (\ref{p2}))), (for the other notations and the terminologies referring to the last section (6)),
 we specialize $(C(\mathbb{F}[S_{d}]),\circ_{q,t})$ to a   commutative associative algebra  $(C(\hat{\mathbb{F}}[S_{d}]),\circ_{A|B})$ with the unit, which will be proven to be isomorphic to the middle-dimensional $\mathbb{\mathbb{C}^*}$-equivalent cohomological rings via the Jack functions over the Hilbert scheme points of $\mathbb{C}^2$ constructed by W. Li, Z. Qin, and W. Wang in \cite{[LQW2]}.}}

\section{Notations and Preliminary}
Throughout this paper, we shall adopt the same notations as I. G.  Macdonald's book \cite{[M]}. Those experts who are familiar with the book can skip this section.

\subsection{Partitions}
A partition \cite{[M]}, p1, is a  (finite or infinite) sequence
\begin{equation}\label{(1.1)}
\lambda=(\lambda_1, \lambda_2,\cdots, \lambda_i, \cdots),
\end{equation}
of non-negative integers in decreasing order:
\begin{equation}
\lambda_1 \geq \lambda_2 \geq \cdots \geq \lambda_i \geq  \cdots \geq 0,
\end{equation}
and containing only finitely many nonzero terms.
The number of the nonzero terms is called the length of $\lambda$, denoted by $l(\lambda)$.  The weight $|\lambda|$ of $\lambda$ is defined by

\begin{equation}
|\lambda|:=\sum_i \lambda_i.
\end{equation}\\

It is   convenient to denote the partition by
\begin{equation}
\lambda= (1^{m_1} 2^{m_2}\cdots r^{m_r}\cdots),
\end{equation}
where $m_i=m_i(\lambda):=\#\{\lambda_j \in \lambda|\lambda_j=i\}$ is called the multiplicity of $i$ in $\lambda.$  Moreover we  define \cite{[M]}, p24
\begin{equation}\label{normal}
z_{\lambda}=\prod_{i\geq 1}i^{m_i}.m_i!
\end{equation}
 If $|\lambda|=d,$ we say that $\lambda$ is a partition of $d,$ and denote it by $ \lambda \vdash d.$
The set of all partitions  of $d$ is
denoted by  ${\mathbb{P}}_d,$   and the set of all partitions (including $d=0$) by ${\mathbb{P}},$   which has a natural lexicographic ordering, \cite{[M]}, p6, i.e.,  for any two partition $\lambda, \mu \in {\mathbb{P}} $ we call  $\lambda \geq \mu$ iff the first non-vanishing difference  $\lambda_i- \mu_i$ is non-negative.

 The Ferrers diagram of a partition $\lambda$ may be formally defined as the set of points
$s=(i, j)\in {\mathbb{Z}}^2$ such that $1 \leq i \leq l(\lambda), 1 \leq j \leq \lambda_i$.
The conjugate of a partition  $\lambda$ is the partition  $\lambda^{\prime}$ whose diagram is the
transpose of the diagram $\lambda$, i.e. the diagram obtained by reflection in the
main diagonal. Hence  $\lambda_i^{\prime}$, is the number of nodes in the $i$-th column of  $\lambda$, or
equivalently
\begin{equation}
\lambda^{\prime}_i = \#\{ j:\lambda_j\geq i\}.
\end{equation}

We will also denote the Ferrers diagram by $\lambda.$
For each $s=(i, j)\in \lambda,$ we define  the arm-length $a(i,j)$, leg-length  $l(i,j)$ \cite{[M]}, p337
and  hook length $hook(i,j)$ \cite{[M]}, p10 by
\begin{eqnarray}
&&a(s)=a(i,j)=\lambda_i-j,\\
&&l(s)=l(i,j)=\lambda^{\prime}_j-i,\\
&&h(s)=h(i,j)=a(i,j)+l(i,j)+1.
\end{eqnarray}

Let $p=(p_1,p_2,p_3,\cdots,)$ be indeterminate variables or the homogeneous symmetric power sums (see, section(\ref{mi})), which are called the
time-variables or Miwa variables, and we assume that $\Delta=(\delta_{1},\cdots,\delta_n)$ is a partition with length $n$. We denote

\begin{eqnarray}
&&\Delta!:=\prod_{i\ge 1}m_i(\Delta)!,\\
&&p_\Delta:=p_{\delta_1}\cdots p_{\delta_n},\label{013}\\
&&\frac{\partial}{\partial p_\Delta}:=\frac{1}{\Delta!}\frac{\partial}{\partial p_{\delta_1}}\cdots \frac{\partial}{ \partial p_{\delta_n}}.\label{partial}
\end{eqnarray}

\subsection {Symmetric function }\label{2.2}
The symmetric group  $S_N$  acts the
 ring ${\mathbb{Q}}[x_1, \cdots, x_N]$ of polynomials in $N$ independent variables
$x=(x_1, \cdots, x_N)$ with rational coefficients by permuting the variables, and a polynomial is symmetric if it
is invariant under this action. The symmetric polynomials form a graded subring:

\begin{eqnarray}
\Lambda_N &=&{\mathbb{Q}}[x_1, \cdots, x_N]^{S_N}\nonumber\\
&=&\oplus_{d\geq 0} \Lambda_N^d,
\end{eqnarray}
where $\Lambda_N^d$  consists of the homogeneous symmetric polynomial of degree $d \geq 0.$

It is well-known that \cite{[M]}, p18, p63, in $\Lambda^d_N,$  the basis consisted by  the monomial symmetric functions $\{m_\lambda(x_1,\cdots,x_N)|\lambda \in {\mathbb{P}}_d\}$ form a  $\mathbb{Z}$-basis, and the element symmetric polynomials  $\{e_{\lambda}(x_1,\cdots,x_N)|\lambda \in {\mathbb{P}}_d\}$ is   dual to the basis consisted by the complete symmetric polynomials $\{h_{\lambda}(x_1,\cdots,x_N)| \lambda \in {\mathbb{P}}_d\},$   the basis consisted by the power sum symmetric polynomials $\{p_{\lambda}(x_1,\cdots,x_N)|\\ \lambda \in {\mathbb{P}}_d\}$ form an orthogonal basis of $\Lambda_N^d$, \cite{[M]}, p64, and the basis consisted by the Schur symmetric polynomials $\{s_{\lambda}(x_1,\cdots,x_N)| \lambda \in {\mathbb{P}}_d\}$ form an orthonormal basis of $\Lambda_N^d$. In the rest of the paper, we will neglect the variables $(x_1,\cdots,x_N)$, and write the above symmetric polynomial as  $m_\lambda, e_\lambda, h_\lambda, p_\lambda, s_\lambda.$

Taking the inverse limit \cite{[M]}, p18,  we have $\mathbb{Z}$-modules $\Lambda^d:$
\begin{equation}
\Lambda^d=\lim_{\leftarrow_N} \Lambda_N^d.
\end{equation}
Let
\begin{equation}
\Lambda=\oplus_{d \geq 0} \Lambda^d,
\end{equation}
so that $\Lambda$ is a free graded $\mathbb{Z}$-module generated by the monomial symmetric polynomial  $m_{\lambda}$ for all partitions $\lambda$. Moreover, the graded ring $\Lambda$ thus defined is called the ring of
symmetric functions in countably many independent variables $x_1, x_2,$ ... . For  any commutative ring or field $\mathcal{A}$, we define
\begin{equation}\label{tensor}
\Lambda_{\mathcal{A}} \cong \Lambda \otimes_{{\mathbb{Z}}} \mathcal{A}.
\end{equation}
Thus by the definition (\ref{013}), $\{p_{\Delta}|\Delta\in {\mathbb{P}}\}$ form a ${\mathbb{Q}}-$
basis of $\Lambda_{\mathbb{Q}}.$

\subsection {Macdonald  functions  with two   parameters $(q,t)$}\label{2.3}
Let $q, t$ be independent variables, and let $\mathbb{ F} = {\mathbb{Q}}(q, t)$ be the field of
rational functions in $q$ and $t$, and let $\Lambda_{\mathbb{ F}} \cong \Lambda \otimes_{{\mathbb{Z}}}{\mathbb{ F}}$ denote the $\mathbb{ F}$-algebra of symmetric functions with coefficients in $\mathbb{ F}$. We define $(q,t)$-deformation of the formula (\ref{normal}) by \cite{[M]}, p309, (2.1)
\begin{equation}\label{43}
 z_{\lambda}(q,t)=z_{\lambda}\prod_{i=1}^{l(\lambda)}\frac{{1-q^{\lambda_i}}}{{1-t^{\lambda_i}}}.
\end{equation}

Then the Macdonald symmetric functions $P_{\lambda}(x; q, t)$
depending rationally on two parameters $(q,t)$ are characterized by the following two properties \cite{[M1]}, \cite{[M]}, p322
 \begin{itemize}
 \item[(a)] There exist $u_{\lambda \mu} \in {\mathbb{F}}$ such that
 \begin{equation}
 P_{\lambda}(x; q, t)=m_{\lambda} + \sum_{\mu  \textless \lambda}u_{\lambda \mu} m_{\mu},
 \end{equation}
 i.e., the transition matrix that expresses the Macdonald functions in
terms of the monomial symmetric functions is a strictly upper unit triangular matrix.
\item[(b)] The  Macdonald functions $P_{\lambda}(x; q, t)$  are pairwise orthogonal relative to the scalar product   defined by
\begin{equation}\label{scalar}
 <p_{\lambda},p_{\mu} >_{q,t}=\delta_{\lambda \mu}{{z}}_{\lambda}(q,t).
\end{equation}
\end{itemize}

\subsection {The  Integral  Macdonald functions  with two parameters $(q,t)$ }
For each partition $\lambda$, we denote \cite{[M]}, p352, (8.1), (8.1$^\prime$),
\begin{eqnarray}
c_\lambda(q, t)=\prod_{s\in \lambda}(1-q^{a(s)}t^{l(s)+1}),\\
c_\lambda^{\prime}(q, t)=\prod_{s\in \lambda}(1-q^{a(s)+1}t^{l(s)}),\\
j_{\lambda}(q, t)=c_\lambda(q, t)c_\lambda^{\prime}(q, t).\label{j}
\end{eqnarray}
taking $q=t^\alpha,$ and letting $ t\to 1$, we define \cite{[M]}, p381, (10.21), (10.23),
\begin{eqnarray}
&&c_\lambda^{(\alpha)}=\lim_{t \to 1}\frac{c_\lambda(t^\alpha, t)}{(1-t)^{|\lambda|}}=\prod_{s\in \lambda}(\alpha a(s)+l(s)+1),\label{path01}\\
&&c_\lambda^{\prime (\alpha)}=\lim_{t \to 1}\frac{c_\lambda^\prime(t^\alpha, t)}{(1-t)^{|\lambda|}}=\prod_{s\in \lambda}(\alpha(a(s)+1)+l(s)),\label{path02}\\
&&j^{(\alpha)}_{\lambda}=\lim_{t \to 1} \frac{j_\lambda(q, t)}{(1-t)^{2|\lambda|}}
=\prod_{s\in \lambda}(\alpha a(s)+l(s)+1)(\alpha(a(s)+1)+l(s)). \label{path1}
\end{eqnarray}

{\remark In general, $j_{\lambda}(q, t)$ is not a perfect square except that $q=t,$ thus we could not use the square root of  $j_{\lambda}(q, t)$ to define the Macdonald-Hurwitz number because we hope that the Macdonald-Hurwitz numbers are the rational functions in the parameters $q,t,$ (see the definition \ref{definition}).

}
{\definition({\cite{[M]}, p352, (8.3))
The integral  Macdonald functions  with two  parameter $(q,t)$ is defined by
\begin{eqnarray}
J_{\lambda}(x; q, t)= c_\lambda(q, t)P_{\lambda}(x; q, t).
\end{eqnarray}}}

Then we have the following orthogonal lemma, which is a key point to our paper,\\
{\lemma {{\bf{the Orthogonal Lemma/Cutting Lemma}} \label{orth} \cite{[M]}, p353, (8.7); p357.  The integral Macdonald functions $J_{\lambda}(x; q, t)$ have scalar products
relative to the scalar product of the formula (\ref{scalar})
\begin{equation}
 <J_{\lambda}(x; q, t),J_{\mu}(x; q, t) >_{q,t}=\delta_{\lambda \mu}j_{\lambda}(q, t).
\end{equation}}}

\remark{{The terminology "cutting" comes from the symplectic surgery theory, for the detail, referring to  \cite{[LZZ]}, \cite{[LR]}.}}

We list  two important special cases of the integral Macdonald functions which are mostly  related to our paper:
\begin{itemize}
\item[(1)] When $ q = t$,
\begin{equation}
J_\lambda(x; t, t)=s_{\lambda}\prod_{(i,j)\in \lambda}(1-t^{hook(i,j)}),
\end{equation}
where $s_{\lambda}$ is the Schur function.

\item[(2)] Taking $q = t^{\alpha}$,  and letting $t \to 1, $ we arrive at  the  Jack symmetric functions, which  were first defined by the statistician Henry Jack in 1969,   \cite{[Jack]}, \cite{[M]}, p353, p381, (10.22),
\begin{equation}\label{path}
J^{(\alpha)}_\lambda(x)=\lim_{t \to 1} \frac{J_\lambda(x; t^{\alpha},t)}{(1-t)^{|\lambda|}}.
\end{equation}
\end{itemize}

{\remark The  Jack functions $J^{(\alpha)}_\lambda(x)$ are the bridge to connect the (integral) Macdonald functions with the zonal functions ($\alpha =2, O(N, \mathbb{R})$), the Schur functions ($\alpha =1, U(N,\mathbb{C})$) and the symplectic zonal functions ($\alpha =\frac{1}{2}, Sp(N, \mathbb{H})$), \cite{[M1]}, \cite{[M]}.}

In our paper, we shall focus on the relationships between the integral Macdonald functions with the power sum functions.
Because the power sum symmetric functions $\{p_{\Delta|\Delta \in {\mathbb{P}}}\}$ form an ${\mathbb{F}}$-basis of $\Lambda_{\mathbb{F}}$, and hence we may express the  integral Macdonald function $J_{\lambda}(x; q, t)$ in terms of them, say
\begin{equation}
J_{\lambda}(x; q, t)=\sum_{\Delta}a_{\lambda}(\Delta)(q,t)p_{\Delta},
\end{equation}
where $a_{\lambda}(\Delta)(q,t) \in {\mathbb{ F}},$ especially, we denote $dim\lambda=a_\lambda(1^{|\lambda|})(q,t),$ which is not generally the unit and very important to construct the Macdonald-Hurwitz numbers.

\remark We believe that $dim\lambda=a_\lambda(1^{|\lambda|})(q,t)$ is some kind of the normalized dimension of the irreducible representation, in fact, for the case $q=t,$ which is actually the normalized dimension of the irreducible complex representation symmetrical group $S_{|\lambda|}$ corresponding to the partition $\lambda.$ For the  case of the Jack functions, $dim\lambda$ is the unit, \cite{[M]}, p382, (10.29).

\bigskip

Then we can restate the Orthogonal Lemma/Cutting Lemma  as \\
{\lemma
\label{oa}
{\begin{itemize}
\item[(1)]
 \bf{The first orthogonal lemma}
\begin{equation}
\qquad \sum_{\{\Delta||\Delta|=|\lambda|\}}a_{\lambda}(\Delta)(q,t) z_{\Delta}(q,t) a_{\mu}(\Delta)(q,t)
=\delta_{\lambda \mu}j_{\lambda}(q, t).
\end{equation}
\item[(2)]
\bf{The second orthogonal lemma}
\begin{equation}\label{second}
\sum_{\{\lambda||\lambda|=|\Delta_1|=|\Delta_2|\}}a_{\lambda}(\Delta_1)(q,t) \frac{1}{j_\lambda(q, t)}a_\lambda(\Delta_2)(q, t)=\delta_{\Delta_1, \Delta_2}\frac{1}{z_{\Delta_1}(q, t)}.
\end{equation}
\end{itemize}
}}

{\remark{ For the orthogonality of the  Jack functions, please refer to  \cite{[M]}, p382, (10.31).}}

\subsection{The initial values}
In this section, we introduce an  initial value of the differential equations (\ref{3}).\\

\bigskip
If $a$ is  indeterminate, we denote by $(a; q)_{\infty}$ the formal  infinite product
\begin{equation}
(a; q)_{\infty} = \prod_{i=0}^{\infty}(1 - aq^i)
\end{equation}
regarded as a formal power series in $a$ and $q$.
Let now $x = (x_1, x_2, ...)$ and $y = (y_l, y_2 ...)$ be two sequences of independent
indeterminate, we define \cite{[M]}, p309, (2.5),
\begin{equation}
\Pi(x,y;q,t)= \prod_{i,j} \frac{(tx_iy_j;q)_\infty}{(x_iy_j;q)_\infty}.
\end{equation}
Then we have \cite{[M]}, p324, (4.13); p309, (2.6)
\begin{eqnarray}
\Pi(x, y; q, t) &=& \sum_{\lambda}  j_\lambda(q, t)^{-1}J_\lambda(x; q,t)J_\lambda(y; q,t), \label{value1}\\
&=&\sum_{\Delta}  z_\Delta(q, t)^{-1}p_\Delta(x)p_\Delta(y), \label{value2}
\end{eqnarray}
where $J_\lambda(x;q,t)$ and $J_\lambda(y; q,t)$ is the integral Macdonald functions in variables $x$
and $y,$ respectively. In the Lemma(\ref{perfect0}), we will give out a geometric proof for the above last two equalities.

\section{Macdonald-Hurwitz number}
After a long preparation, we could define our Macdonald-Hurwitz numbers in this section.\\

Let $\Sigma^h$ be a closed connected Riemann surface of genus $h\geq 0,$ and $\Sigma^g$
 a closed connected Riemann surface of genus $g \geq 0$.
 For a general given points set $\{q_1,\cdots,q_k\}\in \Sigma^h$, which is called the set of
 branch points, we call a
holomorphic map
$f:\Sigma^g\to \Sigma^h$
a ramified
covering of $\Sigma^h$ of degree $d\ge 0$ by $\Sigma^g$ with a ramification type $(\Delta_{1},\cdots,\Delta_{k})$, if
the preimages of $f^{-1}(q_i)=\{c_i^{1}, \cdots, c_i^{l_i}\}$  with orders $\Delta_i = (\delta^i_1,\cdots,\delta^i_{l_i})\vdash d$ for $i=1,\cdots,k,$  respectively, where $l_i=l(\Delta_i)$ is the length of the partition $\Delta_i.$

{\definition\label{definition}
{
We  define the connected genus $g \geq 0$ Macdonald-Hurwitz number  with degree $d > 0$ over the  closed  connected  Riemann surfaces $\Sigma^h$  with the ramification type $(\Delta_{1},\cdots,\Delta_{k})$ in terms of the element of the transformation matrix from the integral Macdonald functions to the power sum  as
\begin{eqnarray}\label{JH}
\hspace{2cm} MH_{h}^{g,d}(\Delta_{1},\cdots,\Delta_{k})=
\sum_{\lambda \vdash d}
(\frac{dim^2\lambda}{j_{\lambda}(q,t)})^{1-h}\frac{a_{\lambda}(\Delta_1)(q,t)}{dim \lambda}\cdots \frac{a_{\lambda}(\Delta_k)(q,t)}{dim \lambda},
\end{eqnarray}
under the condition
\begin{equation}\label{hurwitz}
(2-2h)d-(2-2g)=\sum_{i=1}^{k}(d-l(\Delta_i)),
\end{equation}
otherwise, we define
\begin{eqnarray}{}
MH_{h}^{g,d}(\Delta_{1},\cdots,\Delta_{k})=0,
\end{eqnarray}
and once we take $d=0,$  it is convenient to   define
\begin{eqnarray}
MH_{h}^{g,0}(\Delta_{1},\cdots,\Delta_{k})=0.
\end{eqnarray}
}}

\bigskip

To define the disconnected Macdonald-Hurwitz number, we take the rule of the exponential-logarithmic correspondence, i.e., the disconnected Macdonald-Hurwitz numbers are defined by the product of all connected the Macdonald-Hurwitz numbers module their automorphisms, and then sum all possibilities up, i.e., for any integer $g\in \mathbb{Z},$
    \begin{equation}
    MH^{g,d}_{h}(\Delta_{1},\cdots,\Delta_{k})_{disconnected}=\sum_{m \geq 1}\sum\frac{1}{|Aut|}\prod_{i=1}^m MH^{g_i,d_i}_{h}(\Delta_1^{i_1},\cdots,\Delta_k^{i_k}),
    \end{equation}
    where the sum satisfies the following conditions: for a fixed component number $m$ of the Macdonald-Hurwitz numbers, we have
    \begin{eqnarray}
&&g_i \geq 0, \textit{for $ i=1,\cdots, m$}\\
&&d_i \geq 0,  \textit{for $i=1,\cdots, m$}\\
    &&\Delta_1^{i_1}, \cdots, \Delta_k^{i_k} \vdash d_i \textit{or 0, for $i=1,\cdots, m$}\\
  &&   (2-2h)d-(2-2g_i)=\sum_{j=1}^{k}(d_i-l(\Delta_j^{i_j})),  \textit{for $i=1,\cdots, m$}\\
       && \Delta_j=\cup_{i=1}^m\{\Delta_j^{i_j}\},  \textit{for $j=1,\cdots, k$}\label{sub}\\
                 &&1-g=\sum_{i=1}^m(1-g_i)\\
   && d=\sum_{i=1}^md_i,
     \end{eqnarray}
     where the formula (\ref{sub}) means that  the partitions $\Delta_j^{i_j}$ for  $i=1, \cdots, m,$ compose the partition  $\Delta_j$, i.e., the ramification type of the $j$-th   branch point $q_j$ is invariant when we split the Riemann surface $\Sigma^g$ into m components,  ($\Delta_j^{i_j}$ may be a empty partition, i.e, $\Delta_j^{i_j} \vdash 0,$ the $j$-th branch point $q_j$ might not have pre-image in the $i$-th component.),
    and $|Aut|$ is the cardinality of the automorphisms of all connected components. It is easy to check that the disconnected Macdonald-Hurwitz numbers are the exponential functions of the connected ones.
    To simply our notation, we shall denote the Macdonald-Hurwitz numbers by $ MH^{g,d}_{h}(\Delta_{1},\cdots,\Delta_{k})$ neglecting the connectedness or the disconnectedness.

{\remark{If we take the limit along the path as $q=t^\alpha, t \to 1$, the Macdonald functions will be specialized to the Jack functions, and thus the Macdonald-Hurwitz numbers can be referred to as the Jack-Hurwitz numbers. This is similar to the $b$-Hurwitz numbers ($b=\alpha -1$) discussed in \cite{[CD]}.
}}

By the first orthogonal Lemma,  we can iteratively calculate the Macdonald-Hurwitz numbers by decreasing the genus of the target surface or the numbers of the branch points.

{\theorem\  \label{decrease}
{\begin{itemize}
\item[(1)] Decrease the genus of the target surface directly, i.e., if $h\geq 1,$ we have
\begin{eqnarray}{}
MH_{h}^{g,d}(\Delta_{1},\cdots,\Delta_{k})=\sum_{\Delta \vdash d}MH_{h-1}^{g,d}(\Delta_{1},\cdots,\Delta_{k},\Delta,\Delta)z_{\Delta}(q,t),
\end{eqnarray}
\item[(2)](Associativity) Decrease the genus of the target surface and the numbers of the branched points by cutting it into two parts, i.e, for $h_1+h_2=h, h_1\geq 0, h_2 \geq 0,$ and any positive integer $k> l\ge 1$, we have
\begin{eqnarray}\label{6}
&& MH_{h}^{g,d}(\Delta_{1},\cdots,\Delta_{k})\nonumber\\
&=&
\sum_{\Delta}MH_{h_1}^{g_1,d}(\Delta_{1},\cdots,\Delta_{l},\Delta)z_{\Delta}(q,t)\\
 &&\times MH_{h_2}^{g_2,d}(\Delta,\Delta_{l+1},\cdots,\Delta_{k}).\nonumber
\end{eqnarray}
\end{itemize}
where the genus  $g,g_1,g_2$ are determined by the formula $(\ref{hurwitz})$.}}\\
\begin{proof}
\begin{itemize}
\item[(1)] If $h\geq 1,$ by the first orthogonal lemma,  we have
\begin{eqnarray}{}
&&MH_{h}^{g,d}(\Delta_{1},\cdots,\Delta_{k})\nonumber\\
=&&\sum_{\lambda \vdash d}
(\frac{dim^2\lambda}{j_{\lambda}(q,t)})^{1-h}\frac{a_{\lambda}(\Delta_1)(q,t)}{\dim\lambda}\cdots \frac{a_{\lambda}(\Delta_k)(q,t)}{dim\lambda}\nonumber\\
=&&\sum_{\lambda \vdash d}
(\frac{dim^2\lambda}{j_{\lambda}(q,t)})^{1-h}\frac{a_{\lambda}(\Delta_1)(q,t)}{dim\lambda}\cdots \frac{a_{\lambda}(\Delta_k)(q,t)}{dim\lambda}\nonumber\\
&&\times\sum_{\Delta,\mu}a_{\lambda}(\Delta)(q,t) z_{\Delta}(q,t) a_{\mu}(\Delta)(q,t)
\frac{\delta_{\lambda \mu}}{j_{\lambda}(q, t)}\\\nonumber
=&&\sum_{\lambda \vdash d}
(\frac{dim^2\lambda}{j_{\lambda}(q,t)})^{1-(h-1)}\frac{a_{\lambda}(\Delta_1)(q,t)}{dim\lambda}\cdots \frac{a_{\lambda}(\Delta_k)(q,t)}{dim \lambda}\nonumber\\
&&\times\sum_{\Delta}\frac{a_{\lambda}(\Delta)(q,t)}{dim\lambda} z_{\Delta}(q,t) \frac{a_{\lambda}(\Delta)(q,t)}{dim\lambda}\nonumber\\
=&&\sum_{\Delta \vdash d}MH_{h-1}^{g,d}(\Delta_{1},\cdots,\Delta_{k},\Delta,\Delta)z_{\Delta}(q,t).\nonumber
\end{eqnarray}
\item[(2)]
By the definition of the Macdonald-Hurwitz (\ref{JH}) and the first orthogonal lemma, we have
\begin{eqnarray*}
\mathrm{RHS}&=&\sum_{\Delta}\sum_{\lambda}
(\frac {dim^2\lambda}{j_\lambda (q,t)})^{1-h_1}\frac{a_{\lambda}(\Delta_1)(q, t)}{dim\lambda} \cdots \frac{a_{\lambda}(\Delta_l)(q, t)}{dim\lambda}
\frac{a_{\lambda}(\Delta)(q,t)}{dim\lambda}z_{\Delta}(q,t)\\
&&\times
\sum_{\mu}(\frac {dim^2\mu}{j_\mu(q,t)})^{1-h_2}\frac{a_{\mu}(\Delta)(q, t)}{dim\mu}
\frac{a_{\mu}(\Delta_{l+1})(q, t)}{dim\mu}\cdots \frac{a_\mu(\Delta_k)(q, t)}{dim\mu}\\
&=&\sum_{\lambda}
(\frac {dim^2\lambda}{j_\lambda (q,t)})^{2-h_1-h_2}\frac{a_{\lambda}(\Delta_1)(q, t)}{dim\lambda} \cdots \frac{a_{\lambda}(\Delta_l)(q,t)}{dim\lambda} \delta_{\lambda \mu}\frac{j_\lambda (q,t)}{dim^2\lambda}\\
&&\times
\sum_{\mu}
\frac{a_{\mu}(\Delta_{l+1})(q, t)}{dim\mu}\cdots \frac{a_\mu(\Delta_k)(q, t)}{dim\mu}\\
&=& LHS.
\end{eqnarray*}
\end{itemize}
\end{proof}
{\corollary \label{example0}
{ For any positive integer $d,$ let $(d)\vdash d$ be a partition of $d,$ then we have
\begin{equation}\label{level0}
MH_{0}^{0,d}((d),(d))=\frac{1}{z_{(d)}(q,t)}=\frac{1-t^d}{d(1-q^d)},
\end{equation}
especially,  taking $d=1,$ we have
\begin{equation}\label{level3}
MH_{0}^{0,1}((1),(1))=\frac{1}{z_{(1)}(q,t)}=\frac{1-t}{1-q}.
\end{equation}}}
\begin{proof}
 It follows from the  Theorem \ref{decrease} (2).
\end{proof}

{\remark \label{oo}{
If there is only one branched point in the case of degree one, we have
\begin{equation}\label{level1}
MH_{0}^{0,1}((1))=\frac{1-t}{1-q},
\end{equation}
which follows from that $J_{(1)}(x;q,t)=(1-t)P_{(1)}(x;q,t)=(1-t)p_1$,  \cite{[M]}, p323, (4.8), which does not conflict with the classical Hurwitz number, because in the special case, we have first  to deal $j_{(1)}(q,t)$ and $a_{(1)}((1))(q,t),$ i.e.,
\begin{eqnarray}
j_{(1)}(q,t)=c_{(1)}(q,t)c^\prime_{(1)}(q,t)=(1-t)(1-q),\\
j_{(1)}(1,1)=\lim_{t\to 1}\frac{j_{(1)}(t,t)}{(1-t)^2}=1.
\end{eqnarray}
More generally, if taking $q=t^\alpha,$ and letting $t \to 1,$ the formula (\ref{level1}) would be understood as
\begin{equation}\label{level2}
MH_{0}^{0,1}((1))(\alpha)=\lim_{t \to 1}\frac{1-t}{1-q}=\frac{1}{\alpha},
\end{equation}
where $MH_{0}^{0,1}((1))(\alpha)$ is  the Jack-Hurwitz numbers with the Jack parameter $\alpha.$

{\example \label{example2}
{By the Remark (\ref{oo}), for $d \geq 1,$ $MH_0^{1-d, d}((1^d))=\frac{1}{d!}\frac{(1-t)^d}{(1-q)^d}$, which is a disconnected Macdonald-Hurwitz number (with $m=d$ same components).}}

{\remark{
It is obvious that $a_{\lambda}(\Delta)(q,t)=0$  unless $|\Delta|=|\lambda|$  \cite{[M]}, p356. If  $|\Delta|\textgreater |\lambda|,$ we could take $a_{\lambda}(\Delta)(q,t)=0$, however, if  $|\Delta|\textless|\lambda|,$ we could shift the weight of $\Delta$ by adding $|\lambda|-|\Delta|$-numbers 1, i.e., we get a new partition $\Delta^{\prime}=1^{m_1(\Delta)+(|\lambda|-|\Delta|)}2^{m_2(\Delta)}\cdots$, if $\Delta=1^{m_1(\Delta)}2^{m_2(\Delta)}\cdots$, such that $|\Delta^\prime|=|\lambda|$, then we could define $a_{\lambda}(\Delta)(q,t):={{m_1(\Delta)}  \choose{|\lambda|-|\Delta|+m_1(\Delta)} }a_{\lambda}(\Delta^{\prime})(q,t)$,  \cite{[MMN]},\cite{[Z2]}, then we would obtain a shift Macdonald-Hurwitz number, but in this paper, we shall always assume $|\Delta|=|\lambda|$.}}

{\remark \label{special}\
{\begin{itemize}
\item[(1)] The equation (\ref{hurwitz}) is  exact the Hurwitz formulas.
\item[(2)] Taking   $ q = t$, and letting $t \to 1$, we have
\begin{eqnarray}
&&J^{(1)}_{\lambda}(x) \nonumber\\
&&=\lim_{t \to 1}\frac{J_\lambda(x; t, t)}{(1-t)^{|\lambda|}}\\ \nonumber
&&=\lim_{t \to 1}\frac{s_{\lambda}\prod_{s\in \lambda}(1-t^{hook(s)})}{(1-t)^{|\lambda|}}\\\nonumber
&&=s_{\lambda}\prod_{s\in \lambda}hook(s)\\ \nonumber
&&=\sum_{\Delta}\left[\frac{|\lambda|!}{dim(\lambda)}\frac{\chi_{\lambda(\Delta)}}{z_{\Delta}}\right]p_{\Delta},
\end{eqnarray}
where $s_{\lambda}$ is the Schur function, and $dim(\lambda),$  $\chi_\lambda(\Delta)$ are the dimensions and the character of the irreducible complex representation of symmetric group $S_d$ associated with  $\lambda$ and $\Delta,$ respectively,  \cite{[M]}, \cite{[MMN]}, \cite{[Z1]}. Thus we have
\begin{eqnarray}
&&j_{\lambda}(1,1)=(\frac{|\lambda|!}{dim(\lambda)})^2, \label{key1}\\
&&a_{\lambda}(\Delta)(1,1)=\frac{|\lambda|!}{dim(\lambda)}\frac{\chi_{\lambda}(\Delta)}{z_{\Delta}},\label{key2}\\
&&a_{\lambda}(1^{|\lambda|})(1,1)=\frac{|\lambda|!}{dim(\lambda)}\frac{\chi_{\lambda}(1^{|\lambda|})}{z_{(1^{|\lambda|})}}=1.\label{key2}
\end{eqnarray}
Substitute the above data (\ref{key1})-(\ref{key2}) into (\ref{JH}), we return to the classical Hurwitz number, \cite{[OP]}, \cite{[Z1]}, p4, (8), it is an exact reason to make us choose the integral Macdonald functions with two parameters as our objects, rather than the Macdonald functions.  Moreover, we could obtain that, under the condition (\ref{hurwitz}), the Macdonald-Hurwitz numbers are nonzero, since the classical Hurwitz numbers are nonzero due to the Riemannian existence theorem.
The other reason for us to choose the integral Macdonald functions is that the integral Macdonald functions have better integrality than the Macdonald functions, for example, the normal of the integral Macdonald functions are the polynomial functions in the parameters $q, t, $ however, the normal of the Macdonald functions are the rational functions in the parameters $q,t,$ i.e., \cite{[M]}, p323 (4.11), p339 (6.19), p353. (8.7))
    \begin{eqnarray}
    &&<J_\lambda(x;q,t),J_\lambda(x;q,t)>_{q,t}\\
    =&&c_\lambda(q, t)c^\prime _\lambda(q, t)\nonumber\\
   =&& \prod_{s\in \lambda}(1-q^{a(s)}t^{l(s)+1})(1-q^{a(s)+1}t^{l(s)}),\nonumber\\
   && <P_\lambda(x;q,t),P_\lambda(x;q,t)>_{q,t}\\
   =&&\frac{c^{\prime}_\lambda(q, t)}{c_\lambda(q, t)}\nonumber\\
   =&&\prod_{s\in \lambda}\frac{(1-q^{a(s)+1}t^{l(s)})}{(1-q^{a(s)}t^{l(s)+1})}.\nonumber
    \end{eqnarray}
\item[(3)] Due to the theorem (\ref{decrease})(1), we mainly consider the special case $h=0$ in the rest of the paper, and denote the Macdonald-Hurwitz numbers $MH_{0}^{g,d}(\Delta_{1},\cdots,\Delta_{k})$ by $MH^{g}_d(\Delta_{1},\cdots,\Delta_{k}).$
\end{itemize}}}

\section{The genus-expanded  cut-and-join operators/ integral Macdonald functions with two parameters $(q,t)$}
In this section, we assume all partitions are the partitions of $d.$ Let $\hbar$ be indeterminate, then we have
{\definition \label{opd} {Following the idea of \cite{[Z1]}, for any  partition $\Delta \vdash d$,  we define   a genus-expanded cut-and-join operator $
D(\Delta,\hbar)$ by
\begin{equation}\label{cut}
D(\Delta,\hbar)=\sum_{\Gamma,\Gamma^{\prime}\vdash d}
z_{\Gamma^\prime}(q,t)\hbar^{d+l(\Gamma^{\prime})-l(\Delta)
-l(\Gamma)}MH_d^{g}(\Gamma^{\prime},\Delta,\Gamma)
p_{\Gamma}\frac{\partial}{\partial p_{\Gamma^{\prime}}},
\end{equation}
where  genus $g$ satisfies the  formula (\ref{hurwitz}) ($h=0$):
\begin{equation}
2d-(2-2g)=(d-l(\Gamma^{\prime}))+(d-l(\Delta))+(d-l(\Gamma)),
\end{equation}}}
and $\frac{\partial}{\partial p_{\Gamma^{\prime}}}$ is defined by the formula (\ref{partial}).
\bigskip
{\remark
\begin{itemize}
\item[(1)] Taking the connected Macdonald-Hurwitz number in the above operator (\ref{cut}), and  $\Delta=(1^{d-2}2)$ and $q=t,$ letting $t \to 1,$ we return the classical cut-and-join operator  \cite{[GJ1]}, \cite{[GJ2]}, \cite{[GJ3]}, \cite{[LZZ]}, \cite{[MMN]} or the Laplace-Beltrami operator (\ref{opp}) up a multiple $\frac{1}{2},$ thus the above operator (\ref{cut}) can be thought as a generalization of the Laplace-Beltrami operator (\ref{opp}).
\item[(2)] As in \cite{[Z1]}, the power $d+l(\Gamma^{\prime})-l(\Delta)-l(\Gamma)$ of $\hbar$ in the formula (\ref{cut})
is the "lost" genus after we have cut the target  Riemann surface $\Sigma^h$ into two parts, referring to the Remark(5.8).
\end{itemize}
\bigskip
\definition {Similar to \cite{[Z1]}, we define the genus-expanded integral Macdonald
 functions $J_{\lambda}(x; q, t, \hbar)(\lambda \vdash d)$  by
\begin{equation}\label{11}
J_{\lambda}(x; q, t, \hbar):= \sum_{\Delta \vdash d} \hbar^{-d-l(\Delta)} a_{\lambda}(\Delta)(q,t) p_{\Delta}.
\end{equation}}

{\remark
\begin{itemize}
\item[(1)] If we take $\hbar=1,$ we return to the integral Macdonald functions, thus we have extended the range of the definition of the integral Macdonald functions.
\item[(2)] If we take $q=t,$ and let $t \to 1$, we return to the genus-expanded Schur functions, which first appeared in \cite{[Z1]}.
\end{itemize}
}
\bigskip

 {\theorem \label{eigen} For any partition $\Delta, \lambda \vdash d,$ the following equality holds:
 \begin{equation}\label{12}
 D(\Delta,\hbar)J_\lambda(x; q, t, \hbar) = \hbar^{d-l(\Delta)}a_\lambda(\Delta)(q,t) J_\lambda(x; q, t, \hbar).
\end{equation}
i.e., the genus-expanded integral Macdonald functions  $J_\lambda(x;q,t,\hbar)$
 are  the common  eigenfunctions of  the genus-expanded cut-and-join operators $D(\Delta,\hbar)$
and  $\hbar^{d-l(\Delta)}a_\lambda(\Delta)(q, t)$  are the corresponding eigenvalues.}

\begin{proof} First of all,  for any $\Delta^{\prime}, \Gamma^\prime \vdash d,$ we notice that
\begin{equation}\label{4}
\frac{\partial}{\partial p_\Delta}p_{\Delta^{\prime}}=\delta_{\Delta,\Delta^{\prime}},
\end{equation}
\begin{equation}\label{14}
D(\Delta,\hbar)p_{\Gamma^{\prime}}=\sum_{\Gamma\vdash d}z_{\Gamma^\prime}(q,t)\hbar^{d+l(\Gamma^{\prime})-l(\Delta)-l(\Gamma)}MH_d^{g}(\Gamma^{\prime},\Delta,\Gamma)p_{\Gamma}.
\end{equation}
Then the theorem follows from the first orthogonal lemma and straightforward calculation.
\end{proof}

\bigskip
{\theorem \label{operator}
For a given $d$, as operators acting on  the functions in the time-variables
$p = (p_1, p_2,\cdots )$, the genus-expanded cut-and-join operators are closed under the	multiplication, i.e.,  we have
\begin{equation}\label{19}
 {D(\Delta_1,\hbar) D(\Delta_2,\hbar) = \sum_{\Delta_3}\hbar^{d-l(\Delta_1)-l(\Delta_2)+l(\Delta_3)}
C_{\Delta_1\Delta_2}^{\Delta_3} D(\Delta_3,\hbar) },
\end{equation}
where
\begin{equation}
C_{\Delta_1\Delta_2}^{\Delta_3}=z_{\Delta_3}(q,t)MH_{d}^{g}(\Delta_{1},\Delta_{2},\Delta_{3})
\end{equation}
is a new structure constant
of the center subalgebra $C(\mathbb{F}[S_{d}])$  of the group algebra $\mathbb{F}[S_{d}]$ defined by the formula (\ref{constant}).\\
\begin{proof}
For any $\Gamma^{\prime}\vdash d$, by formula (\ref{14}), we have
\begin{eqnarray}\label{13}
&& \hspace{3cm} D(\Delta_1,\hbar) D(\Delta_2,\hbar)p_{\Gamma^{\prime}} \\
&=&D(\Delta_1,\hbar)\sum_{\Gamma}z_{\Gamma^{\prime}}(q,t)\hbar^{d+l(\Gamma^{\prime})-l(\Delta_2)-l(\Gamma)}MH_{d}^{h_2}(\Gamma^{\prime},\Delta_2,\Gamma)p_{\Gamma}\nonumber\\
&=&\sum_{\Delta}\hbar^{2d+l(\Gamma^{\prime})-l(\Delta_1)-l(\Delta_2)-l(\Delta)}z_{\Gamma^{\prime}}(q,t)\sum_{\Delta^{\prime}}MH_{d}^{g_1}(\Gamma^{\prime},\Delta_2,\Delta^{\prime}){z_{\Delta^{\prime}(q,t)}}MH_{d}^{g_2}(\Delta^{\prime},\Delta_1,\Delta)p_{\Delta} \nonumber\\
&=&\sum_{\Delta}\hbar^{2d+l(\Gamma^{\prime})-l(\Delta_1)-l(\Delta_2)-l(\Delta)}z_{\Gamma^{\prime}}(q,t)MH_{d}^{g_3}(\Gamma^{\prime},\Delta_2,\Delta_1,\Delta)p_{\Delta}\nonumber\\
&=&\sum_{\Delta}\hbar^{2d+l(\Gamma^{\prime})-l(\Delta_1)-l(\Delta_2)-l(\Delta)}z_{\Gamma^{\prime}}(q,t)\sum_{\Delta_3}MH_{d}^{g_4}(\Delta_1,\Delta_2,\Delta_3)z_{\Delta_3}(q,t)MH_{d}^{g_5}(\Delta_3,\Gamma^{\prime},\Delta)p_{\Delta}\nonumber\\
&=&\sum_{\Delta}\hbar^{2d+l(\Gamma^{\prime})-l(\Delta_1)-l(\Delta_2)-l(\Delta)}z_{\Gamma^{\prime}}(q,t)\sum_{\Delta_3}z_{\Delta_3}(q,t)MH_{d}^{g_4}(\Delta_{1},\Delta_{2},\Delta_{3})MH_{d}^{g_5}(\Delta_3,\Gamma^{\prime},\Delta)p_{\Delta}\nonumber\\
&=&\sum_{\Delta_3}\hbar^{d-l(\Delta_1)-l(\Delta_2)+l(\Delta_3)}C_{\Delta_1\Delta_2}^{\Delta_3}D(\Delta_3,\hbar)p_{\Gamma^{\prime}},\nonumber
\end{eqnarray}
which is equivalent to the formula (\ref{19}), and where $g_1, g_2, g_3, g_4,g_5$ are the genuses determined by the formula $(\ref{hurwitz})$.
\end{proof}
\bigskip
{\corollary \label{cor 3.6}
If we normalize the genus-expanded cut-and-join operator $D(\Delta, \hbar)$ by a factor  $\hbar^{-d+l(\Delta)},$
\begin{equation}\label{20}
\hat{D}(\Delta, \hbar):=\hbar^{-d+l(\Delta)}D(\Delta, \hbar),
\end{equation}
then for a given $d$, as operators acting on the space of functions in the time-variables $p=(p_1, p_2, \cdots)$,
all genus-expanded  cut-and-join operators $\hat{D}(\Delta, z)$ for $\Delta \vdash d$ form a commutative associative
algebra, which is denoted by ${{D}}_d,$
 \begin{equation}
 {\hat{D}(\Delta_1,\hbar) \hat{D}(\Delta_2,\hbar) = \sum_{\Delta_3}C_{\Delta_1\Delta_2}^{\Delta_3} \hat{D}(\Delta_3,\hbar), }
\end{equation}\label{CWW0}
 i.e., we have an algebraic isomorphism:
\begin{eqnarray}
{{D}}_d &\cong &(C(\mathbb{F}[S_{d}]),\circ_{q,t})\nonumber\\
        \hat{D}(\Delta,z)   & \mapsto & C_{\Delta},
\end{eqnarray}
where $(C(\mathbb{F}[S_{d}]),\circ_{q,t})$ is the  new commutative associative algebra  defined in theorem (\ref{ca}).
{\begin{proof}
It follows by the straightforward calculation, we skip it.
\end{proof}
}

{\remark It is worth noticing that to deal with the $Schr\ddot{o}dinger$ equation for the one-dimensional N-body problem, the physicist  F. Calogero \cite{[C]} has proposed a similar Laplace-Beltrami operator as the formula (\ref{Delta}) in 1969, it is a little challenge to find out the physics meaning of the genus-expanded cut-and-join operators (\ref{cut}).}

\section {Generating wave functions and its differential equations}
{\definition{ For any given genus $h \geq 0$, degree $d\geq 0$ and partitions $\Delta_1, \cdots, \Delta_n \vdash d,$   the generating wave functions of the degree $d$ is defined by
\begin{eqnarray}\label{99}
&& \hspace{1.5cm} \Phi_h\{\hbar|(u_1,\Delta_1),\cdots,(u_n,\Delta_n)|p^{(1)},\cdots,p^{(k)},p\} \\
&& \qquad =\sum_{l_1,\cdots,l_n\ge 0}\sum_{\Gamma,\Gamma_1,\ldots,\Gamma_{k}\vdash d}\hbar^{2g-2}\nonumber\\
&&\times MH_{h}^{g,d}(\underbrace{\Delta_{1},\cdots,\Delta_1}_{l_1},\cdots,\underbrace{\Delta_{n},\cdots,\Delta_n}_{l_n},\Gamma_1,\cdots,
\Gamma_k,\Gamma)\left[\prod_{j=1}^{n}\frac{u_j^{l_j}}{l_j!}\right]\left[\prod_{i=1}^{k}p^{(i)}_{\Gamma_i}\right]p_\Gamma \nonumber\\
&& \qquad =\sum_{l_1,\cdots,l_n\ge 0}\hbar^{2g-2}(\frac{dim^2\lambda}{j_{\lambda}(q,t)})^{1-h}\left[\prod_{\makebox {j=1}}^n(\frac{a_{\lambda}(\Delta_j)}{dim\lambda})^{l_j}\frac{(u_j)^{l_j}}{l_j!}\right] \left[\prod_{i=1}^{k}\frac{J_{\lambda}(p^{(i)})}{dim\lambda}\right]\frac{J_{\lambda}(p)}{dim\lambda}, \nonumber
\end{eqnarray}
where $\hbar, u_1,\cdots, u_n$ are indeterminate variables,  $p^{(1)}, \cdots, p^{(k)}, p$ are the series of (possibly different) time-variables, and $J_{\lambda}(p^{(i)}), J_{\lambda}(p)$ are the integral Macdonald functions but with the different  series of time-variables, and
$2g-2$ is determined  by the formula (\ref{hurwitz}):
\begin{equation}\label{g1}
(2-2h)d-(2-2g)=\sum_{j=1}^{n}l_j(d-l(\Delta_j))+\sum_{j=1}^{k}(d-l(\Gamma_j))+(d-l(\Gamma)).
\end{equation}
}}

{\remark For the case of the Jack functions, a similar  generating function  appears in \cite{[CD]}, p2, (2).}

{\lemma \label{perfect0}{
Applying the Corollary (\ref{example0}) and the Example (\ref{example2}), we have the following special initial values  (see the equations (\ref{value1})):
\begin{eqnarray}
\Phi_0\{\hbar||p \}&=&\sum_{\lambda}\sum_{\Delta}\hbar^{-d-l(\Delta)}\frac{dim^2\lambda}{j_{\lambda} (q,t)}\frac{a_{\lambda}(\Delta)(q,t)}{dim\lambda}p_{\Delta}\nonumber\\
&=&\sum_{\lambda}\frac{dim^2\lambda}{j_{\lambda} (q,t)}J_{\lambda}(x; q, t, \hbar) \label{algebraic1} \\
&=&\hbar^{-2d}\frac{p_1^d}{d!}\frac{(1-t)^d}{(1-q)^d};\label{geometric1}
\end{eqnarray}

\begin{eqnarray}
\Phi_0\{\hbar||p^{(1)},p\}&=&\sum_{\lambda}\sum_{\Delta_1,\Delta_2}\hbar^{-l(\Delta_1)-l(\Delta_2)}\frac{dim^2\lambda}{j_{\lambda} (q,t)}\frac{a_{\lambda}(\Delta_1)(q,t)}{dim\lambda}\frac{a_{\lambda}(\Delta_2)(q,t)}{dim\lambda}p^{(1)}_{\Delta_1}p_{\Delta_2}\nonumber\\
&=&\sum_{\lambda}\hbar^{2d}\frac{1}{j_{\lambda} (q,t)}J_{\lambda}(q,t,\hbar)(p^{(1)})J_{\lambda}(q,t,\hbar)(p) \label{algebraic}\\
&=&\sum_{\Delta}\hbar^{-2l(\Delta)}\frac{1}{z_{\Delta}(q,t)}p^{(1)}_\Delta p_\Delta \label{geometric},
\end{eqnarray}
where $J_{\lambda}(q,t,\hbar)(p^{(1)})$ and $J_{\lambda}(q,t,\hbar)(p)$ to denote the genus-expanded integral Macdonald function $J_{\lambda}(x;q,t,\hbar),J_{\lambda}(y; q, t, \hbar)$, but with the comprehensible different time-variables series $p^{(1)}, p$.}
\begin{proof}
 The formula (\ref{algebraic1}), (\ref{geometric1}), and (\ref{algebraic}) follow from the definition of Macdonald-Hurwitz numbers and Remark (\ref{oo}). The formula (\ref{geometric}) can be deduced from the Corollary (\ref{example0}) and the restricted condition (\ref{hurwitz}) of the nonzero of the Macdonald-Hurwitz numbers.
We can also prove the above equality directly by the second orthogonal lemma (the formula (\ref{second})) as
\begin{eqnarray}
\Phi_0\{\hbar||p \}&=&\sum_{\lambda}\sum_{\Delta}\hbar^{-d-l(\Delta)}\frac{dim^2\lambda}{j_{\lambda} (q,t)}\frac{a_{\lambda}(\Delta)(q,t)}{dim\lambda}p_{\Delta}\nonumber\\
&=&\sum_{\Delta}\hbar^{-d-l(\Delta)}p_{\Delta}\sum_\lambda\frac{dim^2\lambda}{j_{\lambda} (q,t)}\frac{a_{\lambda}(\Delta)(q,t)}{dim\lambda}\frac{a_\lambda(1^d)}{\dim\lambda}\\
&=&\sum_{\Delta}\hbar^{-d-l(\Delta)}p_{\Delta}\delta_{\Delta, 1^d}\frac{1}{z_{\Delta}(q,t)}\label{conj1}\\
&=&\hbar^{-2d}\frac{p_1^d}{d!}\frac{(1-t)^d}{(1-q)^d};
\end{eqnarray}

\begin{eqnarray}
\Phi_0\{\hbar||p^{(1)},p\}&=&\sum_{\lambda}\sum_{\Delta_1,\Delta_2}\hbar^{-l(\Delta_1)-l(\Delta_2)}\frac{dim^2\lambda}{j_{\lambda} (q,t)}\frac{a_{\lambda}(\Delta_1)(q,t)}{dim\lambda}\frac{a_{\lambda}(\Delta_2)(q,t)}{dim\lambda}p^{(1)}_{\Delta_1}p_{\Delta_2}\nonumber\\
&=&\sum_{\Delta_1,\Delta_2}\hbar^{-l(\Delta_1)-l(\Delta_2)}p^{(1)}_{\Delta_1}p_{\Delta_2}\sum_{\lambda}\frac{1}{j_{\lambda} (q,t)}a_{\lambda}(\Delta_1)(q,t)a_{\lambda}(\Delta_2)(q,t)\\
&=&\sum_{\Delta_1,\Delta_2}\hbar^{-l(\Delta_1)-l(\Delta_2)}p^{(1)}_{\Delta_1}p_{\Delta_2}\delta_{\Delta_1, \Delta_2}\frac{1}{z_\Delta(q,t)}\label{conj3}\\
&=&\sum_{\Delta}\hbar^{-2l(\Delta)}\frac{1}{z_{\Delta}(q,t)}p^{(1)}_\Delta p_\Delta.
\end{eqnarray}
\end{proof}}

{\remark \label{perfect}
\begin{itemize}\
\item[(1)] For the case of the  Jack functions,  taking $q=t^\alpha$, and letting $t \to 1,$  the formula (\ref{geometric1}) and (\ref{geometric}) should be understood as (see also for the hypergeometric functions \cite{[A]}, p6),
     \begin{eqnarray} \label{right}
&&\Phi_0\{\hbar||p\}(\alpha)\nonumber\\
=&&\lim_{t \to 1}\hbar^{-2d}\frac{p_1^d}{d!}\frac{(1-t)^d}{(1-q)^d}\nonumber\\
=&&\hbar^{-2d}\frac{p_1^d}{d!}\frac{1}{{\alpha}^d};
\end{eqnarray}
and
\begin{eqnarray}
&&\Phi_0\{\hbar||p^{(1)},p\}(\alpha)\nonumber\\
=&&\lim_{t \to 1}\sum_{\Delta}\hbar^{-2l(\Delta)}\frac{1}{\frac{z_{\Delta}(q,t)}{(1-t)^d}}p^{(1)}_\Delta p_\Delta \nonumber\\
=&&\sum_{\Delta}\hbar^{-2l(\Delta)}\frac{1}{z_{\Delta}\alpha^{l(\Delta)}}p^{(1)}_\Delta p_\Delta.
\end{eqnarray}

  \item[(2)] If we take $\hbar =1,$ then the calculation from the formula (\ref{algebraic}) to the formula (\ref{geometric}) is equivalent to the calculation from the formula (\ref{value1}) to the formula (\ref{value2}) of the degree  $d,$ thus we give out a geometric proof for the calculation from the formula (\ref{value1}) to the formula (\ref{value2}).
\item[(3)] The formula (\ref{geometric1}) can be regarded  as the genus $g$, the degree $d$ and parameters $q,t$-expansions of hook formula \cite{[M]}, \cite{[MMN]},  and the formula (\ref{geometric}) as  the genus $g$, the degree $d$ and parameter $q,t$-expansions  of the refined Cauchy-Littlewood identity \cite{[BW]},\cite{[M]},\cite{[MMN]}, \cite{[MMN2]}.
\end{itemize}
}

\example
{Let us check the formula ((\ref{geometric1}), (\ref{geometric}) for the case of  the  integral Macdonald functions of the degree 2,  we have
\begin{eqnarray}
&&J_{(2)}(x; q,t)=\frac{1}{2}(1+q)(1-t)^2 p_{1^2} + \frac{1}{2}(1-q)(1-t^2)p_2;\\
&&J_{(1^2)}(x; q,t)=\frac{1}{2}(1-t)(1-t^2) p_{1^2} - \frac{1}{2}(1-t)(1-t^2)p_2;\\
&&j_{(2)}(q,t)=(1-t)(1-qt)(1-q)(1-q^2);\\
&&j_{(1^2)}(q, t)=(1-t)(1-t^2)(1-qt)(1-q);
\end{eqnarray}
Thus, it is easy to check that the formula (\ref{geometric1}), (\ref{geometric}) is correct.}

\example
{Let us check the formula (\ref{geometric1}), (\ref{geometric}) for the case of  the  Jack functions of the degree 2,  we have
\begin{eqnarray}
&&J_{(2)}^{(\alpha)}(x)=p_1^2 + \alpha p_2;\\
&&J_{(1^2)}^{(\alpha)}(x)=p_1^2 - p_2;\\
&&j_{(2)}^{(\alpha)}=2\alpha^2(\alpha + 1);\\
&&j_{(1^2)}^{(\alpha)}=2\alpha(\alpha + 1).
\end{eqnarray}
Thus, it is easy to check that the formula(\ref{geometric1}), (\ref{geometric})  is correct.}

\bigskip
{\theorem \label{de}For any $i$, we have
\begin{eqnarray}\label{3}
&&\nonumber \frac{\partial\Phi_h\{\hbar|(u_1,\Delta_1),\cdots,(u_n,\Delta_n)|p^{(1)},\cdots,p^{(k)},p\}}{\partial u_i} \nonumber\\
& =&D(\Delta_i,\hbar)\Phi_h\{\hbar|(u_1,\Delta_1),\cdots,(u_n,\Delta_n)|p^{(1)},\cdots,p^{(k)},p\}.
\end{eqnarray}}
\begin{proof}
Firstly, we have
$$
\frac{\partial\Phi_g\{\hbar|(u_1,\Delta_1),\cdots,(u_n,\Delta_n)|p^{(1)},\cdots,p^{(k)},p\}}{\partial u_i}$$
$$=\sum_{l_1,\cdots,l_n\ge 0}
\sum_{\Gamma_1,\ldots,\Gamma_{k}, \Gamma}\hbar^{2g-2}MH_{h}^{g,d}(\underbrace{\Delta_{1},\cdots,\Delta_1}_{l_1},\cdots,\underbrace{\Delta_{n},\cdots,\Delta_n}_{l_n},\Gamma_1,\cdots,
\Gamma_k,\Gamma)$$
$$\frac{(u_i)^{l_i-1}}{(l_i-1)!}\left[\prod_{\makebox {j=1,}\\ j\ne i}^n\frac{(u_j)^{l_j}}{l_j!}\right]
\left[\prod_{j=1}^{k}p^{(j)}_{\Gamma_j}\right]p_\Gamma.
$$
We can write RHS of equation (\ref{3}) as
\begin{eqnarray*}
\mathrm{RHS}&=&\sum_{l_i\ge 1} \sum_{l_1,\cdots,\check{l_i},\cdots,l_n\ge 0}
\sum_{\Gamma_1,\ldots,\Gamma_{k},\Gamma^{\prime}}\hbar^{2g^i_--2}\\
&&\times
MH_{h}^{g^i_-,d}(\underbrace{\Delta_{1},\cdots,\Delta_1}_{l_1},\cdots,\underbrace{\Delta_{i},
\cdots,\Delta_i}_{l_i-1},\cdots,\underbrace{\Delta_{n},\cdots,\Delta_n}_{l_n},\Gamma_1,\cdots,
\Gamma_k,\Gamma^{\prime})\\
&&\times
\frac{(u_i)^{l_i-1}}{(l_i-1)!}[\prod_{\makebox {j=1,} j\ne i}^n\frac{(u_j)^{l_j}}{l_j!}]
[\prod_{j=1}^{k}p^{(j)}_{\Gamma_j}]D(\Delta_i,\hbar)p_{\Gamma^{\prime}},
\end{eqnarray*}
where $\check{l_i}$ means that we omit $l_i$, and  $2g^i_--2$ is also determined by the formula (\ref{hurwitz}):
\begin{equation}
(2-2h)d-(2-2g^i_-)=\sum_{j=1}^{n}l_j(d-l(\Delta_j))-(d-l(\Delta_i))+\sum_{j=1}^{k}(d-l(\Gamma_j))+(d-l(\Gamma^{\prime})).\label{g2}
\end{equation}
Moreover we have the following facts:
\begin{itemize}
\item[Fact (1)]: By Theorem (\ref{eigen}),
$$D(\Delta_i,\hbar)p_{\Gamma^{\prime}}=z_{\Gamma^{\prime}}(q,t)\sum_{\Gamma}\hbar^{d+l(\Gamma^{\prime})-l(\Delta_i)-l(\Gamma)}MH_{d}^{g^i_+}(\Gamma^{\prime},\Delta_i,\Gamma)p_{\Gamma};$$
\item[Fact (2)]:
 $MH_{d}^{g^i_+}(\Gamma^{\prime},\Delta_i,\Gamma)\ne 0$ only if
\begin{equation}\label{g3}
2d-(2-2g^i_+)=(d-l(\Gamma^{\prime}))+(d-l(\Delta_i))+(d-l(\Gamma));
\end{equation}
\item[Fact (3)]: By the formula (\ref{g1}), (\ref{g2}), (\ref{g3}), we have
$$g=g^i_++g^i_-+l(\Gamma^{\prime})-1.$$
\end{itemize}
Then the theorem follows from the first orthogonal lemma (\ref{oa}).
\end{proof}

\bigskip
\remark
\begin{itemize}
\item[(1)]
We notice that the  same  phenomenon of "genus lost" appears as \cite{[Z1]}, \cite{[Z2]}, i.e.,
we "lost" the genus after we have executed the  cutting surgery:
\begin{eqnarray}
&&(2g-2)-(2g^i_--2)\\
&=&(2g^i_+-2)+2l(\Gamma^{\prime})\nonumber\\
&=&d+l(\Gamma^{\prime})-l(\Delta_i)-l(\Gamma).\nonumber
\end{eqnarray}
\item[(2)] It is because the generating function $\Phi_h$ satisfies the equation (\ref{3}), we would like to call it by the generating wave function of the degree $d.$ It is obvious that we should add all degree $d$ and all genus $g$ to the infinity to get the classical wave functions, but we would have to deal with the countable infinite variables,  and we would like to deal with the case in the future.
\end{itemize}

 \bigskip
{\corollary {\label{expression}
\begin{eqnarray}
&&\Phi_h\{\hbar|(u_1,\Delta_1)\cdots,(u_n,\Delta_n)|p^{(1)},\cdots,p^{(k)},p\} \nonumber\\
=&&\left[\prod_{i=1}^n exp (u_iD(\Delta_i,\hbar))\right]\Phi_h\{\hbar||p^{(1)},\cdots,p^{(k)},p\}.
\end{eqnarray}
In particular, if we take $h=0,k=0, 1$, we have
\begin{eqnarray}
&&\Phi_0\{\hbar|(u_1,\Delta_1),\cdots,(u_n,\Delta_n)|p\}\nonumber\\
=&&[\prod_{i=1}^{n}exp(u_iD(\Delta_i,\hbar))](\hbar^{-2d}\frac{p_1^d}{d!}\frac{(1-t)^d}{(1-q)^d});
\end{eqnarray}
\begin{eqnarray}
&&\Phi_0\{\hbar|(u_1,\Delta_1),\cdots,(u_n,\Delta_n)|q,p\}\nonumber\\
=&&[\prod_{i=1}^{n}exp(u_iD(\Delta_i,\hbar))](\sum_{\Delta}\hbar^{-2l(\Delta)}\frac{1}{z_{\Delta}(q,t)}q_\Delta p_\Delta).
\end{eqnarray}
}}
\begin{proof}
 It follows from Theorem (\ref{de}) and the initial values (\ref{geometric1}), (\ref{geometric}).
\end{proof}

\section{A new  commutative associative structure over the cental subalgebra $C(\mathbb{F}[S_{d}])$ and its application}

In this section, we use Macdonald-Hurwitz numbers to define a $\circ_{q,t}$-product over the classical central subalgebra $C(\mathbb{F}[S_d])$ of the group algebra $\mathbb{F}[S_d]$,   where $\mathbb{F}=\mathbb{Q}(q,t)$ is the field of rational functions with rational coefficients in the parameters $q,t,$ and $S_d$ is the symmetric group of the degree $d.$ This product can be thought of as a $(q,t)$-deformation of the product in $C(\mathbb{C}[S_d])$ that extends  the coefficient field $\mathbb{C}$ to the field $\mathbb{F}$.

{\theorem \label{ca} { For any partition $\Delta \vdash d, $ let $C_{\Delta} \in C(\mathbb{C}[S_{d}])$ be the classical central element of $C(\mathbb{C}[S_{d}])$ corresponding to $\Delta$, \cite{[LZZ]}, \cite{[Z1]}, then  $\{C_{\Delta}| \Delta \in{\mathbb{P}}_d\}$ also form a ${\mathbb{F}}$-basis of $C(\mathbb{F}[S_{d}])$. By linearly extending,  $C(\mathbb{F}[S_{d}])$ is a new commutative associative algebra with the unit $C_{(1^d)}$ under the following product,
\begin{equation}
C_{\Delta_1}\circ_{q,t}C_{\Delta_2}= C_{\Delta_1\Delta_2}^{\Delta_3}C_{\Delta_3},
\end{equation}}
where
\begin{equation}\label{constant}
C_{\Delta_1\Delta_2}^{\Delta_3}=z_{\Delta_3}(q,t)MH_{d}^{g}(\Delta_{1},\Delta_{2},\Delta_{3})
\end{equation}
is the new structure constant.}
\begin{proof}
It follows from the commutativity (Definition (\ref{JH})) and the associativity (Theorem (\ref{decrease}) (2))  for the case $h=0$) of the Macdonald-Hurwitz numbers.
\end{proof}

{\lemma {{\bf{and Definition}}
Let us introduce a bilinear form $<-,->_{q,t}$  and trilinear form $T_{q,t}(-,-,-)$ over $C(\mathbb{F}[S_{d}])$ as
\begin{eqnarray}
&&<C_{\Delta_1},C_{\Delta_2}>_{q,t}:=\delta_{\Delta_1,\Delta_2}\frac{1}{z_{\Delta_1}(q,t)},\label{adopt1}\\
&&T_{q,t}(C_{\Delta_1},C_{\Delta_2},C_{\Delta_3}):=MH_d^g(\Delta_1,\Delta_2, \Delta_3).\label{adopt2}
\end{eqnarray}}}
Then we have
\begin{eqnarray}
&&T_{q,t}(C_{\Delta_1},C_{\Delta_2},C_{\Delta_3})\\
&&=<C_{\Delta_1},C_{\Delta_2}\circ_{q,t}C_{\Delta_3}>_{q,t}\\
&&=<C_{\Delta_1}\circ_{q,t}C_{\Delta_2}, C_{\Delta_3}>_{q,t}
\end{eqnarray}}}
{\corollary\label{6.1}{ $C(\mathbb{F}[S_{d}])$ is a Frobenius algebra (referring to \cite{[LP]}) with the unit $C_{(1^d)}$ for the bilinear form $<-,->_{q,t}$ and trilinear form $T_{q,t}(-,-,-)$.
Moreover, for $\lambda \vdash d,$
\begin{eqnarray}
\epsilon_\lambda(q,t):=\sum_{\Delta \vdash d} \frac{\dim\lambda}{j_\lambda(q,t)}a_\lambda(\Delta)(q,t)z_\Delta(q,t)C_\Delta
\end{eqnarray}
form the  idempotent bases of $(C(\mathbb{F}[S_{d}]), \circ_{q,t})$, which are orthogonal to each other,
\begin{eqnarray}
<\epsilon_\lambda(q,t),\epsilon_\mu(q,t)>_{q,t}=\delta_{\lambda,\mu}\frac{dim^2\lambda}{j_\lambda(q,t)}.
\end{eqnarray}}}
\begin{proof}
It follows from the definition of the Frobenius algebra and  $\circ_{q,t}$-product and direct calculation.
\end{proof}

{\remark{
Taking  $q=t,$ and letting $t \to 1,$ we return to the classical central subalgebra  $C(\mathbb{C}[S_{d}])$, (It is worth noticing  that the bilinear form  is different from the  standard bilinear form (with normal $z_\Delta(q,t)$, rather than $\frac{1}{z_\Delta(q,t)}$),  thus we give out an $(q,t)$-deformation product of the central  subalgebra $C(\mathbb{C}[S_{d}])$ when we extend the coefficient field $\mathbb{C}$ of $C(\mathbb{C}[S_{d}])$ to the field $\mathbb{F}.$}
}}

  In the rest of this section, we establish a connection between the new commutative and associative algebra $(C(\mathbb{F}[S_d]), \circ_{q,t})$ and the middle-dimensional $\mathbb{C}^*$-equivalent cohomological ring via the Jack function over the Hilbert scheme points of $\mathbb{C}^2$ constructed by W. Li, Z. Qin, and W. Wang \cite{[LQW2]}.

Let $\alpha$ be the Jack parameter of the Jack functions $J^{(\alpha)}_\lambda(x)$, and $\hat{\mathbb{F}}:=\mathbb{Q}(A, B)$ be the field of the rational functions with the rational coefficient in  $A, B,$  where $A, B$ are two indeterminates such that $\alpha=\frac{B}{A}$.

 Let $C(\hat{\mathbb{F}}[S_d]), \circ_{A|B})$ be the restriction of $(C(\mathbb{F}[S_{d}]), \circ_{q,t})$ to the special  path $\eta(A|B)$ defined by
\begin{eqnarray}
&&t=r^A,\label{p1}\\
&&q=r^B,\label{p2}
\end{eqnarray}
where $r$ is the path parameter, and letting $r \to 1$, which is different from the path $q=t^\alpha, t=t \to 1$ in the Macdonald's book \cite{[M]} or the formula (\ref{path}), although we will eventually obtain the same  Jack functions with the Jack parameters $\alpha=\frac{B}{A}$ from the  Macdonald functions in both case. We denote the corresponding quantities by adding $A|B$, for example, we have the following data:

\begin{eqnarray}
&&c_\lambda(A|B):=\lim_{r\to 1}\frac{c_\lambda(r^B, r^A)}{(1-r)^{|\lambda|}}=\prod_{(i,j)\in \lambda}(A(l(i,j)+1)+Ba(i,j)),\label{path03}\\
&&c^\prime_\lambda(A|B):=\lim_{r\to 1}\frac{c^\prime_\lambda(r^B, r^A)}{(1-r)^{|\lambda|}}=\prod_{(i,j)\in \lambda}(Al(i,j)+B(a(i,j)+1)),\label{path04}\\
&&j_{\lambda}(A|B):=\lim_{r\to 1}\frac{j_\lambda(r^B, r^A)}{(1-r)^{2|\lambda|}}\label{path2}=c_\lambda(A|B)c_\lambda^{\prime}(A|B).\label{path05}
\end{eqnarray}
We should notice that the formula (\ref{path03})-(\ref{path05}) is different from the formula (\ref{path01})-(\ref{path1}) up to a factor $A^{|\lambda|}$  or $A^{2|\lambda|}$ respectively,  since we have taken the limit along the different paths. Moreover the above descriptions are the exact means of  "we take the limit of some qualities along a path", or "we take the restrict of a ring to a path".

Moreover we have the bilinear form $<-,->_{A|B}$  and trilinear form $T_{A|B}(-,-,-)$ over $C(\hat{\mathbb{F}}[S_{d}])$ as
\begin{eqnarray}
&&<C_{\Delta_1},C_{\Delta_2}>_{A|B}=\delta_{\Delta_1,\Delta_2}\frac{1}{z_{\Delta_1}\alpha^{l(\Delta_1)}};\\
&&T_{A|B}(C_{\Delta_1},C_{\Delta_2},C_{\Delta_3})\\
&&:=<C_{\Delta_1},C_{\Delta_2}\circ_{A|B}C_{\Delta_3}>_{A|B}\\
&&=<C_{\Delta_1}\circ_{A|B}C_{\Delta_2}, C_{\Delta_3}>_{A|B}\\
&&=MH_d^g(\Delta_1,\Delta_2, \Delta_3)(A|B),
\end{eqnarray}
where $MH_d^g(\Delta_1,\Delta_2, \Delta_3)(A|B)$ is the limit of the Macdonald-Hurwitz numbers $MH_d^g(\Delta_1,\Delta_2, \Delta_3)$ along the path $\eta(A|B)$.

We restate the Corollary (\ref{6.1}) as \\
{\lemma \label{62}{ $C(\hat{\mathbb{F}}[S_{d}])$ is a Frobenius algebra  with the unit $C_{(1^d)}$ for the bilinear form $<-,->_{A|B}$ and trilinear form $T_{A|B}(-,-,-)$.
Moreover, for $\lambda \vdash d,$
\begin{eqnarray}
\epsilon_\lambda(A|B):=\lim_{r \to 1}\epsilon_\lambda(r^B,r^A)(1-r)^{|\lambda|}
\end{eqnarray}
are the  idempotent bases of $(C(\hat{\mathbb{F}}[S_{d}]), \circ_{A|B})$, which are orthogonal to each other,
\begin{eqnarray}
<\epsilon_\lambda(A|B),\epsilon_\mu(A|B)>_{A|B}=\delta_{\lambda,\mu}\frac{1}{j_\lambda(A|B)}.
\end{eqnarray}}}
\begin{proof}
It follows from the Corollary(\ref{6.1}) and $dim \lambda=1$ for the Jack functions,  \cite{[M]}, p382, (10.29).
\end{proof}

In the rest of the paper, we will adopt the same notations as  \cite{[LQW1]}, \cite{[LQW2]}.

Let $\theta$ be the 1-dimensional standard $\mathbb{T}:=\mathbb{C}^\star$-module. Fix two nonzero integers $A$ and $B$  with the same signs, which are corresponding to the parameters $\alpha, \beta$ in \cite{[LQW1]}, \cite{[LQW2]}. Let $u, v$ be
the standard coordinate functions on $\mathbb{C}^2$. For  $s\in \mathbb{T}$, we define the action of $\mathbb{T}$ on $\mathbb{C}^2$ by \cite{[LQW2]}, Example 2.1,
\begin{equation}
s\cdot(u, v) = (s^Au, s^{-B}v).
\end{equation}
The origin of $\mathbb{C}^2$ is the only fixed point, which is denoted by $ O$.

Let $H^\star_{\mathbb{T}}(\mathbb{C}^2)$ be the equivariant cohomology ring of $\mathbb{C}^2$
with $\mathbb{C}$ coefficients.
$H^\star_{\mathbb{T}}(\mathbb{C}^2)$ is a $\mathbb{C}[\bar{t}]$-module if we identify $H^\star_{\mathbb{T}}(pt)$ with $\mathbb{C}[\bar{t}]$ ($\bar{t}$
is an element of the degree -2).

Let $F$ and
$F^\prime$ be the $u$-axis and $v$-axis respectively in $\mathbb{C}^2$. As $\mathbb{T}$-modules, we have $T_OF = \theta^{-A}$ and  $T_OF^\prime = \theta^{B}$,  then we have \cite{[LQW2]}
\begin{eqnarray}
[F]=-A^{-1}\bar{t}^{-1}[O], [F^\prime]=B^{-1}\bar{t}^{-1}[O],
\end{eqnarray}

Let $(\mathbb{C}^2)^{[d]}$ be the Hilbert scheme parameterizing
all the 0-dimensional closed subschemes $\xi$ of $\mathbb{C}^2$ with $dim_{\mathbb{C}} H^0(\mathcal{O}_\xi) = d$. The $\mathbb{T}$-
action on $\mathbb{C}^2$ induces a $\mathbb{T}$-action on  $(\mathbb{C}^2)^{[d]}$. The support of a $\mathbb{T}$-fixed point in  $(\mathbb{C}^2)^{[d]}$
is contained in $(\mathbb{C}^2)^{\mathbb{T}}=O$ and indexed by partitions $\lambda$ of $d$.

Let $\xi_\lambda$ be the fixed points
in $((\mathbb{C}^2)^{[d]})^{\mathbb{T}}$ corresponding to a partition $\lambda$ of $d$. Then for $\lambda \vdash d$, the tangent space
of $(\mathbb{C}^2)^{[d]}$ at the fixed point $\xi_\lambda$ is $\mathbb{T}$-equivariantly isomorphic to \cite{[LQW2]}, p4, (2.10):
\begin{equation}
T_{\xi_{\lambda}}(\mathbb{C}^2)^{[d]}=\oplus_{(i,j)\in \lambda} \theta^{A(l(i,j)+1)+Ba(i,j)} \oplus \theta^{-Al(i,j)-B(a(i,j)+1)}.
\end{equation}

We have the equivariant Euler class \cite{[LQW2]}, p4, (2.11),
\begin{equation}
e_{\mathbb{T}}(T_{\xi_{\lambda}}(\mathbb{C}^2)^{[d]}) =(-1)^d c_\lambda(A|B)c^\prime_\lambda(A|B)\bar{t}^{2d}.
\end{equation}

Note that $\xi_{\lambda} \in H^{4d}_{\mathbb{T}}((\mathbb{C}^2)^{[d]})$.  Then  for $\lambda \vdash d,$   the following distinguished class \cite{[LQW2]}, p4, (2.14):
\begin{equation}
[\lambda] =(-1)^d\frac{1}{c_\lambda(A|B)}\bar{t}^{-d}[\xi_\lambda]
\end{equation}
form a linear basis of the $\mathbb{C}$-vector space $H^{2d}_{\mathbb{T}}((\mathbb{C}^2)^{[d]})$, which have $\star$-product structure \cite{[LQW2]}, p8,
\begin{equation}
\frac{[\lambda]}{c^\prime_\lambda(A|B)} \star \frac{[\mu]}{c^\prime_\mu(A|B)}=\delta_{\lambda,\mu}\frac{[\lambda]}{c^\prime_\lambda(A|B)}, \lambda, \mu \vdash d,
\end{equation}
for $[\lambda], [\mu] \in H^{2d}_{\mathbb{T}}((\mathbb{C}^2)^{[d]})$ and the bilinear form $<-,->$ \cite{[LQW2]}, p8,
\begin{equation}
<[\lambda],[\mu]>=\delta_{\lambda,\mu}\frac{c^\prime_\lambda(A|B)}{c_\lambda(A|B)}.
\end{equation}

 We have the following theorem,
{\theorem{There exists an isomorphism that preserves the bilinear form  between $(C(\hat{\mathbb{F}}[S_{d}]), \circ_{A|B})$ and $(H^{2d}_{\mathbb{T}}((\mathbb{C}^2)^{[d]}),\star)$.
}}
\begin{proof} Applying  the lemma (\ref{62}), let us construct the isomorphism which maps the  idempotent element to the  idempotent element, i.e., for $\lambda \vdash d,$
\begin{eqnarray}
\rho:&&(C(\hat{\mathbb{F}}[S_{d}]), \circ_{A|B})\to (H^{2d}_{\mathbb{T}}((\mathbb{C}^2)^{[d]}),\star)\\
&&\epsilon_\lambda(A|B) \mapsto\frac{[\lambda]}{c^\prime_\lambda(A|B)}.
\end{eqnarray}
Moreover, we have
\begin{eqnarray}
&&<\rho(\epsilon_\lambda(A|B)),\rho(\epsilon_\mu(A|B))>\\
&&=<\frac{[\lambda]}{c^\prime_\lambda(A|B)},\frac{[\mu]}{c^\prime_\mu(A|B)}>\\
&&=\delta_{\lambda,\mu}\frac{1}{j_\lambda(A|B)}\\
&&=<\epsilon_\lambda(A|B),\epsilon_\mu(A|B)>_{A|B}.
\end{eqnarray}
Thus $\rho$ is the  isomorphism we need.
\end{proof}

\remark{{ The equivariant cohomological class $[\lambda]$ can be corresponded to the Jack polynomial,  as detailed in Theorem 3.2 of \cite{[LQW2]}.}}

{\remark{
In our views,  $(H^{2d}_{\mathbb{T}}(\mathbb{C}^2)^{[d]},\star)$  constructed by W. Li, Z. Qin,  and W. Wang should be called $\mathbb{T}$-equivalent cohomological ring via the special path $q=r^B, t=r^A,  r \to 1$ over the Hilbert scheme points of $\mathbb{C}^2$}} (with two parameters $A, B$),  rather than via the Jack functions ( which should take the limit along the path $q=t^\alpha, t \to 1, $ with only one parameter $\alpha,$ in fact,  W. Li, Z. Qin, and W. Wang actually took the limit along the path $\eta(A|B)$ in \cite{[LQW2]}), whose idempotent elements in both cases are different up to a  factor  $A^{|\lambda|}$.

\bigskip
{\bf Acknowledgements} The author would like to thank Prof. Xiaojun Chen, Prof. An-min Li, Prof. Wei Luo and Prof. Yongbin Ruan, Prof. Qi Zhang,  Guanghui Fu,  and Yin Liu for their communications and helps.  This work was partially supported by  NSFC 11890663 and NSFC 12171351, which are led by  Prof. Bohui Chen and Prof. Song  Yang, respectively.

\end{document}